\documentclass[twocolumn,letterpaper]{IEEEAerospaceCLS}  

\usepackage[]{graphicx}    
\newcommand{\ignore}[1]{}  
\usepackage{cite}
\usepackage{amsmath}
\usepackage{amssymb}
\usepackage{comment}
\usepackage{algorithm}
\usepackage{algpseudocode}
\usepackage{nicefrac, xfrac}

\begin{document}
\title{Robust Closed-form Framework for Drag-Propulsive Control of Formation Flight}

\author{%
Matthew Hunter\\ 
Department of Aeronautics \& Astronautics\\
Stanford University\\
Stanford, CA 94305\\
mhunter8@stanford.edu
\and 
Simone D'Amico\\
Department of Aeronautics \& Astronautics\\
Stanford University\\
Stanford, CA 94305\\
damicos@stanford.edu
\thanks{\footnotesize 979-8-3503-0462-6/24/$\$31.00$ \copyright2024 IEEE}              
}

\maketitle

\thispagestyle{plain}
\pagestyle{plain}

\maketitle

\thispagestyle{plain}
\pagestyle{plain}

\begin{abstract}
    Distributed Space Systems (DSS) can accomplish many objectives that would be difficult or impossible with a single monolithic spacecraft, but often use spacecraft with limited propulsion and computing power. To conserve propellant, DSS spacecraft can leverage environmental forces to produce low cost relative forces, such as atmospheric drag in Low Earth Orbit (LEO). Atmospheric drag can be controlled on-board spacecraft by modulating differential attitude to create differential drag with little to no propulsion. This work maximizes this strategy by developing robust, closed-form algorithms for optimal "hybrid" control, defined here as complementing the magnitude and directionality deficiencies of differential drag with the superior control authority of propulsive maneuvers within the same control window. The use of drag introduces a significant source of dynamic uncertainty to the control problem, primarily from inaccurate atmospheric density modelling. This paper develops a robust control architecture to reject dynamic modeling errors and recover final state accuracy performance. The architecture periodically re-solves the remaining drag maneuver plan in a diminishing horizon policy, using the current state estimate to reject operational and navigation errors while retaining the cost savings offered by differential drag. The approach is built upon reachable set theory, a control analysis technique that provides conditions for optimality for Linear Time Variant (LTV) systems with norm-like cost functions and optimally decouples hybrid control into separate propulsive and differential drag sub-problems. The optimal decoupling allows the architecture to solve control for hybrid or drag-only maneuvering. The optimality conditions are leveraged to prove the asymptotic Lyapunov stability of ideal reachable set-based control solutions and derive error bounds on this stability. A covariance analysis on the primary sources of error in hybrid control determines that the dominating source of uncertainty is the atmospheric density model. This conclusion allows the bounds on Lyapunov stability to be completely met by simply conditioning the atmospheric density value used in the control solutions. The derived Lyapunov stability bounds and performance of the algorithm are demonstrated experimentally in high fidelity full-force orbital simulations of the upcoming Space Weather Reconfigurable Multiscale Experiment (SWARM-EX) mission. The overall approach used by this work is not specific to DSS and can be extrapolated to any LTV control system with norm-like cost.
\end{abstract}

\tableofcontents

\section{Introduction}

Formation flight, part of the broader field of Distributed Space Systems (DSS), enables multiple spacecraft to accomplish objectives that would be infeasible for a single monolithic spacecraft, with additional built-in redundancies to single component failure. The GRACE, PRISMA, TanDEM-X, Starling, and MMS missions, among others, have demonstrated the capabilities of distributed instrumentation and navigation in orbit and drive interest in future missions such as the VIrtual Super Optics Reconfigurable Swarm (VISORS) and the Space Weather Atmospheric Reconfigurable Multiscale EXperiment (SWARM-EX), a major motivator of this work \cite{grace,damico_phd,starling,burch_MMS,koenig_VISORS,swarm-ex}. Satellite formations and swarms commonly incorporate smaller, cheaper spacecraft to lower mission cost, introducing restrictions on propellant and collision safety that limit mission lifetime and autonomy. In Low Earth Orbit (LEO), atmospheric drag can supplement propulsion through on-board modulation of differential attitude to generate differential drag at little to no propulsive cost. However, drag generally suffers from poor characterization due to inaccurate atmospheric density modelling. This uncertainty, in addition to that produced by real world propulsion, attitude determination and control, and state estimation, can cause a significant deterioration on the final state accuracy of control solutions and drives interest in balancing cost-efficient maneuvering with robust control strategies that reject operational sources of uncertainty and provide control accuracy guarantees.

Error rejection and robust control techniques are fundamental to the problem of in-flight mission implementation of ideal maneuver planning. Extensive literature exists on the subject, focusing mostly on propulsion given its use on many previous DSS missions, and different feedback techniques are used depending on the formulation of the control problem and the choice of state representation. Propulsive maneuvers may be approximated as impulsive, or instantaneous changes in spacecraft velocity, to create Linear Time-Variant (LTV) models of controlled relative spacecraft motion and kinematically optimal control profiles, for both cartesian and Integration Constant state representations. The CanX-4/-5 mission used a cartesian relative state and Linear Quadratic Regulator (LQR) framework to balance minimum delta-v reconfiguration cost with desired state tracking errors, re-solving the remaining maneuver plan immediately before conducting each maneuver \cite{roth_phd}. Relative Orbital Elements (ROE), an Integration Constant parameterization, offer increased insight over cartesian states into the secular, long-term relative motion of spacecraft by removing the short-term dynamics from the LTV model, resulting in fully deterministic solutions. The PRISMA and TanDEM-X missions used closed-form, delta-v efficient algorithms for regularly scheduled formation keeping, accomplishing error rejection in a similar strategy to CanX-4/-5 by re-solving maneuvers until the moment of actuation \cite{damico_phd}. Lippe et al. and Holzinger et al. applied Lyapunov stability theorems to ROE-based control algorithms to provide guarantees for stable, delta-v efficient maneuvering in near-circular or eccentric orbit and bounds on stable, uncoordinated formation keeping and collision safety \cite{lippe,holzinger}. Chernick and D'Amico also leveraged the ROE LTV formulation to conduct a covariance sensitivity analysis and quantify the effects of dynamic modelling error sources on closed-form optimal maneuver plans \cite{chernick}.

Simple, robust feedback control laws for continuous thrust, with gains defined as a function of current spacecraft state, can be derived from the Algebraic Ricatti Equation solution to LQR, Lyapunov theory, sliding mode control, or nonlinear switching lanes, for both cartesian and ROE states \cite{alfriend,silvestrini,steindorf,cho,koenig_drag}. Given the overall higher complexity of the continuous dynamics model, several works have applied Model Predictive Control (MPC) to incorporate additional guidance-based considerations into the control problem, including a multi-objective cost function, thrust magnitude and directional constraints, and collision safety \cite{nair,breger,lim}. Many of the techniques used for low-thrust propulsive control have been extended to differential drag control, an inherently low-thrust maneuver, to converge despite poor drag characterization. Mazal et al. and Riano-Rios et al. have implemented LQR with provable Lyapunov stability to converge to a desired state with guaranteed error bounds for known uncertainty and with online parameter estimation on the primary sources of uncertainty \cite{mazal,riano-rios_ICL}. Dell'Elce et al. and Hu et al. instead used differential drag in a finite-horizon MPC architecture to accommodate arbitrary constraints while minimizing control window duration or distribute formation regularization throughout a swarm while ensuring reconfiguration constraint feasibility \cite{dellelce,hu}. Koenig et al. demonstrated that a nonlinear state space switching lane control law for low-thrust control was equally applicable to differential drag, rejecting atmospheric density modeling errors by defining a bang-bang and deadband control law only dependent on the spacecraft's current state \cite{koenig_drag}.

A very limited number of works have expanded beyond a single control source to consider the more complex roto-translational problem, producing relative accelerations from both on-board thrusters and differential attitude. The dual-quaternion parameterization combines translational and rotational motion into a single dynamics model, and applications of LQR, sliding mode, and Proportional Integral control have produced asymptotically stable, but suboptimal, controllers \cite{stanfield,zhang,li}. Riano-Rios et al. also accomplished coupled attitude stabilization and changes in translational motion under a single Lyapunov stable controller, but only for a particular spacecraft structural design \cite{riano-rios_roto}. Shouman et al. instead employed a specialized satellite construction to decouple translational and attitude control and demonstrated the ability of differential drag to reject perturbations on propulsive maneuver plans with Lyapunov stability \cite{shouman}. However, none of these roto-translational robust control approaches address optimal maneuvering, specifically using differential drag not only to reject disturbances to the propulsive maneuver plan, but also to decrease the propulsive reconfiguration cost.

To meet this gap, this paper presents a novel DSS "hybrid" maneuver planning architecture, complementing the magnitude and directional deficiencies of differential drag with the superior control authority of propulsive maneuvers within the same control window. The architecture uses the current state estimate to periodically re-solve the remaining maneuver plan and reject maneuver execution, estimation, and dynamic propagation errors while retaining differential drag cost savings. The foundation of this architecture is a generic reachable set theory-based optimal hybrid control solver. Reachable set theory is a control analysis technique implemented in previous work by the authors that provides optimality conditions for both individual maneuvers and the overall maneuver plan, enables the creation of provably-optimal control solvers, and largely decouples hybrid control into mutually optimal drag and propulsion sub-problems \cite{koenig_RST,chernick,hunter_hybrid}. This paper provides three major contributions to the state of the art. First, the re-solve architecture defines an algorithmic structure for stable, convergent, and provably-optimal hybrid DSS control. Second, generic reachable set theory-based provably-optimal solvers are proven to be asymptotically Lyapunov stable, and bounds on maneuver plan error are derived that guarantee asymptotic convergence. Third, the main sources of state error in hybrid control are detailed qualitatively and compared quantitatively through covariance error ellipsoid performance bounds. 

This novel robust hybrid control approach is defined in four major sections. First, the re-solve architecture is established for hybrid and drag-only control. The hybrid setup extends the reachable set-based maneuver planning to a shrinking horizon policy, re-solving drag maneuvers for the remaining control window. The drag-only setup adjusts the control window length to accommodate the limitations on drag control authority and converge to the final state without propellant. Second, the architecture is proven to be asymptotically Lyapunov stable under ideal maneuvering, and bounds on maneuver plan error are derived that guarantee this convergence. Third, covariance error analysis of hybrid control identifies the primary causes of final state performance loss to inform the design of the re-solve architecture, such that all operational errors can be rejected with drag re-solves alone. Finally, the architecture's performance is validated in a high fidelity full-force orbital simulation that propagates the Gauss Variational Equations (GVE) for the relevant LEO perturbations.

\section{Problem Definition}

\subsection{State Representation and Dynamics}

Here, formations and swarms consist of two or more spacecraft flying in close proximity, such that the relative motion between them can be linearized about a reference orbit, traditionally of an uncontrolled spacecraft. This orbit, known as the "chief", can be designated within the formation or can be imaginary, and the motion of the other spacecraft, each designated as "deputy," is defined relative to this orbit. Two state representations are of primary interest in this paper to describe and enable DSS control. The first is the cartesian Radial-Tangential-Normal (RTN) frame, defined from the chief outward radially from the central body, normally along the chief's angular momentum, and tangentially to complete the right-handed triad. The RTN frame is typically used to express deputy control actions. The second and equivalent representation is quasi-nonsingular ROE, a nonlinear combination of the Keplerian orbital elements, defined as 
\begin{equation} \label{eq: ROE}
    \delta\boldsymbol{\alpha} = 
    \begin{bmatrix}
    \delta a \\ \delta \lambda \\ \delta e_x \\ \delta e_y \\ \delta i_x \\ \delta i_y
    \end{bmatrix} = 
    \begin{bmatrix}
    \Delta a/a_c \\
    \Delta u + \Delta \Omega \cos i_c \\
    e_d\cos\omega_d - e_c\cos\omega_c \\
    e_d\sin\omega_d - e_c\sin\omega_c \\
    \Delta i \\
    \Delta \Omega \sin i_c
    \end{bmatrix},
\end{equation} 
\noindent where $\Delta$ is a difference in the associated quantity between the deputy (subscript $d$) and the chief (subscript $c$), $\delta a$ is the relative semi-major axis, $\delta \lambda$ is the relative mean longitude, $\delta e_x$ and $\delta e_y$ are the components of the relative eccentricity vector $\delta \boldsymbol{e}$, $\delta i_x$ and $\delta i_y$ are the components of the relative inclination vector $\delta \boldsymbol{i}$, $u = M + \omega$ is the mean argument of latitude, and $a$, $e$, $i$, $\Omega$, $\omega$, and $M$ are the Keplerian orbital elements. The quasi-nonsingular formulation is valid for orbits of arbitrary eccentricity and in particular for near-circular orbits. This work uses mean ROE, calculated from mean Keplerian orbital elements, to remove short period oscillations from the state. Here, mean refers to the orbit average of J2-perturbed two-body motion. The mean ROE under unperturbed Keplerian dynamics are constant values equivalent to the Integration Constants of the Hill-Clohessy-Wiltshire (HCW) and Yamanaka-Ankersen (YA) differential equations of spacecraft relative motion at small separations \cite{alfriend}. As a function of the Keplerian orbital elements, both osculating and mean ROE vary slowly under orbital perturbations, enabling the creation of an LTV dynamic constraint model of perturbed Keplerian motion with RTN control actions. The secular and long-term evolution of the ROE is captured by a State Transition Matrix (STM) that linearly propagates the state over a desired time interval.  LEO orbits have negligible eccentricity with the primary perturbations of J2 and atmospheric drag. Given that differences in atmospheric drag will be posed as a control force, this work uses an STM that includes corrections for J2 for mean ROE in near-circular orbit, given as
\begin{align} \label{eq: Gamma matrix}
    &\boldsymbol{\Phi}^{J_2}(t) = \boldsymbol{\Phi}^{J_2} (\boldsymbol{\alpha}_c(t),t_f-t) \\
    & \! = \! 
    \begin{bmatrix}
    1 & 0 & 0 & 0 & 0 & 0 \\
    \Phi_{21}^{J_2}(t) & 1 & 0 & 0 & \Phi_{25}^{J_2}(t) & 0 \\
    0 & 0 & \cos(\Dot{\omega}_c\tau) & -\sin(\Dot{\omega}_c\tau) & 0 & 0 \\
    0 & 0 & \sin(\Dot{\omega}_c\tau) & \cos(\Dot{\omega}_c\tau) & 0 & 0 \\
    0 & 0 & 0 & 0 & 1 & 0 \\
    \Phi_{61}^{J_2}(t) & 0 & 0 & 0 & \Phi_{65}^{J_2}(t) & 1 \\
    \end{bmatrix}, \nonumber
\end{align}
\noindent with the following simplifying substitutions
\begin{align}
    &\tau = t_f - t \text{ , } \eta_c = \sqrt{1-e_{c}^{2}} \text{ , } \\
    &\kappa_c = \frac{3}{4}\frac{J_2 R^{2}_{E} \sqrt{\mu_E}}{a_{c}^{7/2} \eta_{c}^{4}} \text{ , } \Dot{\omega}_c = \kappa_c(5\cos^2(i_c)-1) \nonumber \\
    &\Phi_{21}^{J_2}(t) = -(\frac{3}{2}n_c+\frac{7}{2}\kappa_c(1+\eta_c)(3\cos^2(i_c)-1))\tau \nonumber \\
    &\Phi_{25}^{J_2}(t) = -\kappa_c(4+3\eta_c)\sin(2i_c)\tau \text{ , } \nonumber \\
    &\Phi_{61}^{J_2}(t) = \frac{7}{2}\kappa_c\sin(2i_c)\tau \text{ , } \Phi_{65}^{J_2}(t) = 2\kappa_c\sin^2(i_c)\tau, \nonumber
\end{align}
\noindent where $\boldsymbol{\alpha}_c(t)$ are the mean Keplerian orbital elements at time $t$, $t_f$ is the end of the control window, $n$ is the mean motion, $R_E$, $J_2$, and $\mu_E$ are the radius, oblateness, and standard gravitational parameter of Earth, and $\Phi^{J_2}_{i,j} (t)$ is the element in the $i$th row of the $j$th column \cite{koenig_STM}. The map between RTN delta-v control actions and resulting changes in osculating ROE is derived through the Gauss Variational Equations and captured in control input matrix $\boldsymbol{B}$ for near-circular orbit, given as
\begin{equation} \label{eq: B matrix}
    \hspace{-2mm}
    \boldsymbol{B}(t) = 
    \frac{1}{a_c n_c}
    \begin{bmatrix}
    0 \! \! & 2 \! \! & 0 \\
    -2 \! \! & 0 \! \! & 0 \\
    \sin u_c \! \! & 2\cos u_c \! \! & 0 \\
    -\cos u_c \! \! & 2\sin u_c \! \! & 0 \\
    0 \! \! & 0 \! \! & \cos u_c \\
    0 \! \! & 0 \! \! & \sin u_c \\
    \end{bmatrix},
\end{equation}
\noindent where $\boldsymbol{B}(t) = \boldsymbol{B}(\boldsymbol{\alpha}_c(t))$ and $n$ is the mean motion \cite{damico_phd}. It is important to note that $\boldsymbol{B}$ also applies approximately to mean ROE and decouples in-plane (IP) control (radial and tangential maneuvers with $\delta a$, $\delta \lambda$, and $\delta \boldsymbol{e}$ states) from out-of-plane (OOP) control (normal maneuvers with $\delta \boldsymbol{i}$ states) as a product of linearization \cite{roscoe}. The dynamic constraint model of formation flight can be posed as
\begin{equation} \label{eq: dynamics}
    \delta \boldsymbol{\alpha}_f = \boldsymbol{\Phi}^{J_2} (t_0) \delta \boldsymbol{\alpha}_0 + \int_{t_0}^{t_f} \boldsymbol{\Phi}^{J_2} (t) \boldsymbol{B} (t) \boldsymbol{u}(t) dt,
\end{equation}
\noindent with initial state $\delta \boldsymbol{\alpha}_0$ at $t_0$, final desired state $\delta \boldsymbol{\alpha}_f$ at $t_f$, and the RTN actions $\boldsymbol{u}(t)$ occurring in the reconfiguration control window $[t_0,t_f]$. $\boldsymbol{u}$ in this paper will include both propulsive thrust $\boldsymbol{u}^{\text{prop}}$ and differential drag $\boldsymbol{u}^{\text{drag}}$.

\subsection{Differential Drag}

Atmospheric drag, on a single spacecraft, is a non-conservative force that acts directly against the direction of motion, given in RTN as 
\begin{equation} \label{eq: atmospheric drag}
    \boldsymbol{p}^{{\text{drag}}} = 
    \begin{bmatrix}
        0 & -\frac{1}{2}\rho v^2 \frac{C_D A}{m} & 0
    \end{bmatrix}^{\text{T}} = 
    \begin{bmatrix}
        0 & -\frac{1}{2}\rho n^2 a^2 B & 0
    \end{bmatrix}^{\text{T}},
\end{equation}
\noindent where $\rho$ is the atmospheric density, $v = na$ approximates spacecraft tangential velocity in near-circular orbit, $C_D$ is the spacecraft's coefficient of drag, $A$ is the cross-sectional area, $m$ is the spacecraft's mass, and $B = \frac{C_D A}{m}$ defines the spacecraft's ballistic coefficient. The near-circular assumption aligns the spacecraft's velocity tangentially with its RTN frame, simplifying drag to a negative tangential-only force. Spacecraft capable of attitude control can modulate this force by changing $A$ normal to the spacecraft's velocity vector, and differential drag, between a chief and deputy, can be controlled through differential attitude. Differential drag can then be posed as a control force in RTN as
\begin{align} \label{eq: drag force}
    \hspace{-2mm}
    \boldsymbol{u}^{\text{drag}}(t) = 
    \boldsymbol{p}^{\text{drag}}_d(t) - \boldsymbol{p}^{\text{drag}}_c(t) = 
    \begin{bmatrix}
    0 & \! \! \frac{1}{2}n_c^2a_c^2\Delta B^{\rho}(t) & \! \! 0
    \end{bmatrix}^{\text{T}},
\end{align}
\noindent where $\Delta B^{\rho}(t) = \Bar{\rho}_c (B_c(t) - B_d(t))$ defines the augmented differential ballistic coefficient. Atmospheric density can fluctuate greatly over short periods of time, so the orbital time-average $\Bar{\rho}$ is used for simplicity. This assumption has little effect on the output solution given the low instantaneous magnitude of differential drag, but it is not required to find a valid drag maneuver plan. By assuming close orbits and small separations with respect to orbital radius, both the chief and deputy share the chief's $\Bar{\rho}$, $n$, $a$, and RTN frame, since RTN orientation varies by less than $1^\circ$ for along-track separations up to 100km. The magnitude of $\boldsymbol{u}^{\text{drag}}$ is bounded by the limits on $\Delta B^{\rho}(t)$, dictated by the minimum and maximum $A$ for both the chief and deputy. It is worth noting that differential drag is almost always present between any two spacecraft because differential attitude is rarely zero, and residual differences in $B$ exist even for identically built spacecraft \cite{grace}. The effect of differential attitude on spacecraft relative dynamics can be quantified by substituting (\ref{eq: drag force}) into (\ref{eq: dynamics}).

\subsection{Optimal Control Problem}

The general optimal control problem used by this paper is given as
\begin{align} \label{eq: general control problem}
    \text{minimize }& \int_{t_0}^{t_f} f(\boldsymbol{u}(t),t) \\
    \text{subject to }& \delta \boldsymbol{\alpha}_f = \boldsymbol{\Phi}^{J_2} (t_0) \delta \boldsymbol{\alpha}_0 + \int_{t_0}^{t_f} \boldsymbol{\Phi}^{J_2} (t) \boldsymbol{B} (t) \boldsymbol{u}(t) dt, \nonumber
\end{align}
\noindent where $f$ is a norm-like cost function. The attitude constraints of each spacecraft, roto-translational coupling, and deviations from the drag maneuver plan to conduct propulsive maneuvers are all neglected, implying that propulsive and differential drag control can be used independently and simultaneously. Several aspects of this general formulation can be made more specific. First, propulsive actions are assumed to be "impulsive," or instantaneous additions of delta-v, given that burn time is short with respect to orbital period. Second, attitude maneuvering is assumed to require negligible propulsive-cost compared to translational propulsive maneuvering, such that drag is approximated as a zero-cost maneuver in the cost function. Third, the spacecraft is assumed to have a single propulsive thruster, and the appropriate choice of $f$ is an L2-norm of each propulsive maneuver. Fourth, the dynamic constraint is rearranged to group the invariants of the reconfiguration problem into a quantity known as the pseudostate $\Delta \delta \boldsymbol{\alpha}$, given as
\begin{equation} \label{eq: pseudostate}
    \Delta \delta \boldsymbol{\alpha} = \delta \boldsymbol{\alpha}_f - \boldsymbol{\Phi}^{J_2} (t_0) \delta \boldsymbol{\alpha}_0.
\end{equation}
\noindent Finally, the control window is discretized into an arbitrary series of time steps, and differential drag is assumed to be constant over each time step $[t_{j},t_{j+1}]$, such that its effect over each time step can be analytically integrated to pseudostate space. The optimal control problem considered by the remainder of this paper is given as
\begin{equation} \label{eq: hybrid control problem}
    \hspace{-5mm} \text{minimize} \hspace{-0.5mm} \sum_{i=1}^{k} \left\| \boldsymbol{C}^{\text{prop}} \boldsymbol{u}_i \right\|_{2} \hspace{-0.5mm}, \text{subject to } \hspace{-0.5mm} \Delta \delta \boldsymbol{\alpha} \hspace{-0.5mm} = \hspace{-1mm} \sum_{i=1}^{k} \hspace{-0.5mm} \boldsymbol{\Gamma}(t_{i}) \boldsymbol{u}_{i}, \hspace{-4.5mm}
\end{equation}
\noindent for $t_i \in [t_0,t_f]$. Control inputs $\boldsymbol{u}_i = [\boldsymbol{u}_{i}^{\text{propulsive}};\boldsymbol{u}_{i}^{\text{drag}}]$ concatenate both forms of control in the RTN frame, and propulsive cost is isolated by selector matrix $\boldsymbol{C}^{\text{prop}} = [\boldsymbol{I}_3,\boldsymbol{0}_3]$. The auxiliary matrix $\boldsymbol{\Gamma}(t_i) \hspace{-0.5mm} = \hspace{-0.5mm} [\boldsymbol{\Phi}^{J_2}(t_i)\boldsymbol{B}(t_i),\int_{t_{i}}^{t_{i+1}}\boldsymbol{\Phi}^{J_2}(t)\boldsymbol{B}(t)dt]$ translates propulsive maneuvers at time $t_i$ and analytically integrates differential drag maneuvers held constant over time step $[t_{i},t_{i+1}]$ to pseudostate space. The control window $[t_0,t_f]$ is discretized into $k$ time steps, such that $\sum_{i=1}^{k} (t_{i+1} - t_{i}) = t_f - t_0$. A primary advantage of this formulation is the norm-like cost function with LTV constraints, properties that will be leveraged in the next section to derive conditions for control optimality through reachable set theory.

\section{Closed-form Optimal Planning}

\subsection{Reachable Set Theory}

The proposed robust architecture is built upon the conditions for optimality found through reachable set theory, and the overall approach is summarized here. Reachable set theory primarily considers two sets under the dynamics in (\ref{eq: dynamics}): the set $S(c,T)$ of all $\Delta \delta \boldsymbol{\alpha}$ that can be reached by a single maneuver $u$ of non-negative arbitrary cost $c$ over times in the control window $T$ and the set $S^*(c,T)$ of all $\Delta \delta \boldsymbol{\alpha}$ that can be reached by $k \geq 1$ maneuvers $\boldsymbol{u}_i$ at times $t_i \in T$ of combined cost no greater than $c$, given as
\begin{align}
    S(c,T) = \{\Delta \delta \boldsymbol{\alpha}:& \Delta \delta \boldsymbol{\alpha} = \boldsymbol{\Gamma}(t) \boldsymbol{u}, \left\|\boldsymbol{u}\right\|_{2} \leq c, t \in T \}, \label{eq: RST S} \\
    S^*(c,T) = \{\Delta \delta \boldsymbol{\alpha}:& \Delta \delta \boldsymbol{\alpha} = \Sigma^{k}_{i=1} \boldsymbol{\Gamma}(t_i)\boldsymbol{u}_i, t_i \in T, \label{eq: RST S*} \\
    & \hspace{20mm} \left\|\boldsymbol{u}_i \right\|_{2} \leq c_i, \Sigma^{k}_{i=1} c_i = c \}. \nonumber
\end{align}
It can be seen that $S^*$ is a linear combination, and therefore convex hull, of the $\Delta \delta \boldsymbol{\alpha}$ in $S$. Given a norm-like cost function that scales linearly with magnitude, the scalar $c$ that scales $S^*$ until its boundary contains the reconfiguration $\Delta \delta \boldsymbol{\alpha}$ is the minimum or optimal cost $\delta v_{\text{min}}$.

Several observations about this formulation establish the authors' approach to hybrid control \cite{chernick}. The 2$n$-dimensional (2$n$D) reconfiguration can be broken down into $n$ 2D planes, where the 2D plane that contains the highest $\delta v_{\text{min}}$, known as the dominant plane, drives the minimum cost of the full reconfiguration. Within each 2D plane, the reachable set can be broken down into distinct contours known as dominance case candidates, defined by the outward normal of the supporting hyperplane for that contour \cite{koenig_RST}. The dominance case of a reconfiguration is the distinct contour of the reachable set that intersects $\Delta \delta \boldsymbol{\alpha}$ at minimum cost within the dominant plane. The minimum cost $\delta v_{\text{min}}$ for the reconfiguration is found as the solution to
\begin{equation} \label{eq: delta-v min}
    \hspace{-2mm}
    \text{maximize: } \frac{\boldsymbol{\eta}^{\text{T}} \Delta \delta \boldsymbol{\alpha}}{\max_{t \in T} \left( \max_{||\boldsymbol{u}|| \leq 1 } \left( \boldsymbol{\eta}^{\text{T}} \boldsymbol{\Phi}^{J_2}(t) \boldsymbol{B}(t) \boldsymbol{u} \right) \right) },
\end{equation}
\noindent where $\boldsymbol{\eta}$ is the supporting hyperplane of the reconfiguration dominance case \cite{koenig_RST}. Assuming a finite number of dominance case candidates exist, the maximization problem can be solved by quantifying (\ref{eq: delta-v min}) for each dominance case candidate and corresponding $\boldsymbol{\eta}$, the largest of which is the $\delta v_{\text{min}}$ for the reconfiguration. Furthermore, optimal maneuver times $T_{\text{opt}}$ and directions for each dominance case candidate, resulting from the denominator of (\ref{eq: delta-v min}), occur at points on $S^*$, defined by $\boldsymbol{\eta}$, that share a boundary with $S$. 

\begin{figure}
    \centering
    \includegraphics[scale=0.4]{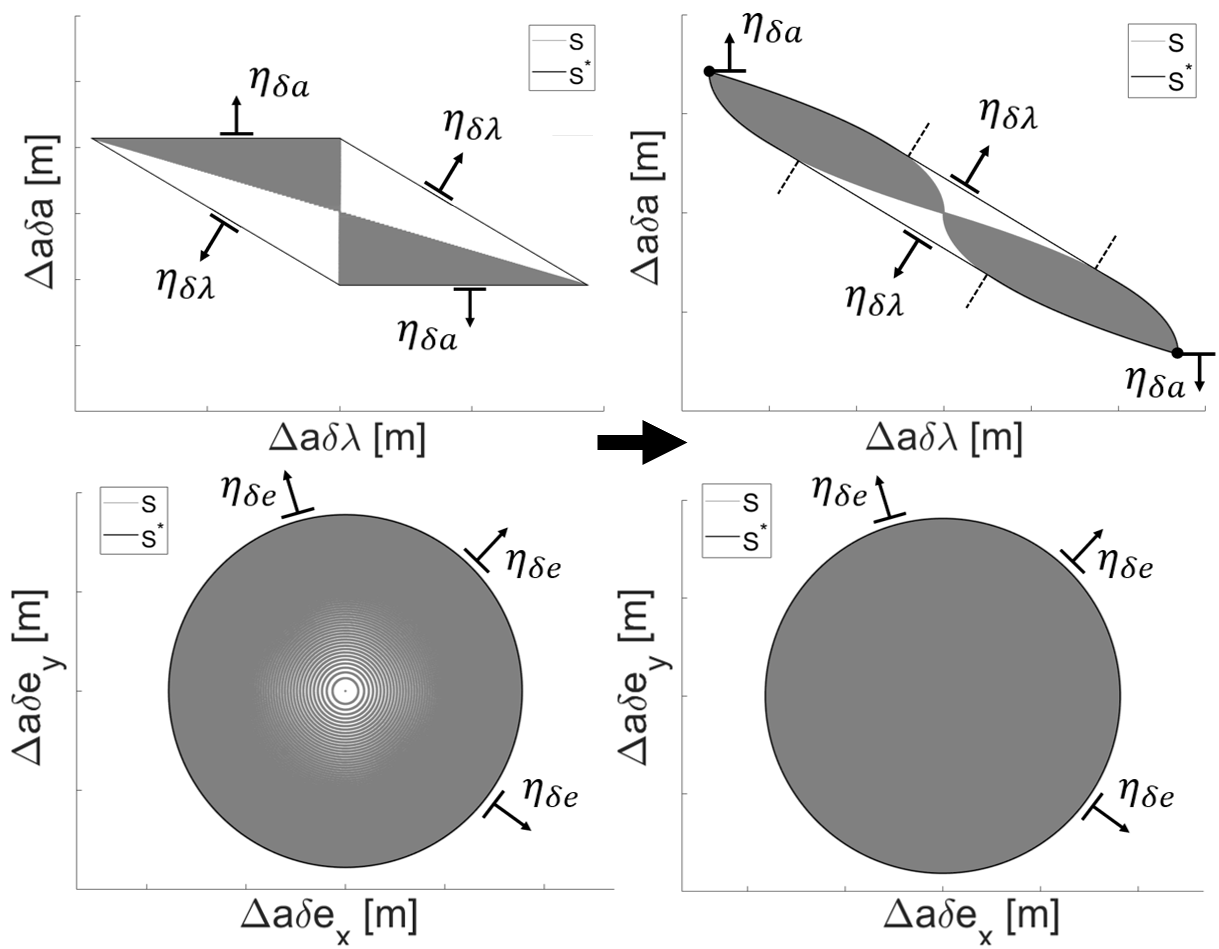}
    \caption{In-plane reachable sets for propulsive control (left column) and differential drag (right column) over a 100 orbit control window in the $\Delta \delta \boldsymbol{a}$ plane (top row) and $\Delta \delta \boldsymbol{e}$ plane (bottom row). Dominance case candidates and related $\boldsymbol{\eta}$ correspond to unique contours of the propulsive $S^*$, and these $\boldsymbol{\eta}$ are applied to the differential drag $S^*$ to identify optimal attitude profiles.}
    \label{fig: IP RS}
\end{figure}

Applying this methodology to the hybrid control problem in (\ref{eq: hybrid control problem}), the hybrid reachable sets are found by taking a single maneuver of arbitrary cost and numerically sampling the $\Delta \delta \boldsymbol{\alpha}$ effect over all possible maneuver times and directions in the control window to find the set $S$ from (\ref{eq: RST S}). For differential drag specifically, the zero-cost maneuvers are sampled over maneuver start times, maneuver durations limited by $t_f$, and possible magnitudes, from the negative to positive bound on $\Delta B^\rho$. The ROE can be conveniently arranged into 3 2D planes to match the decoupling found in $\boldsymbol{B}$: the $\Delta \delta \boldsymbol{a}$ plane ($\Delta \delta a$, $\Delta \delta \lambda$), the $\Delta \delta \boldsymbol{e}$ plane, and the $\Delta \delta \boldsymbol{i}$ plane. This arrangement enables the IP and OOP reconfiguration sub-problems to be analyzed and solved independently, while still accounting for the J2 coupling present in $\boldsymbol{\Phi}^{J_2}(t)$, and the removal of drag from the OOP sub-problem. Propulsive and differential drag control can also be decoupled into independent IP sub-problems. Because drag is a zero-cost maneuver, the dominance case candidates and corresponding $\boldsymbol{\eta}$ that result from the propulsive analysis can also be aligned with the differential drag reachable set to identify optimal drag profiles, or the modulation of differential attitude over time. Fig. \ref{fig: IP RS} illustrates this approach visually for the IP reconfiguration.

\noindent Chernick and D'Amico leveraged the simplicity of the propulsive ROE reachable sets to identify three IP dominance case candidates ($\delta a$, $\delta \lambda$, and $\delta e$) and one OOP dominance case candidate ($\delta i$), named for the dimension that drives minimum cost in the dominance case candidate's contour \cite{chernick}. All IP dominance case candidates have tangential-only optimal maneuvers, allowing analytical expressions of the $\boldsymbol{\eta}$ and $\delta v_{\text{min}}$ for each dominance case candidate. It is worth noting that the $\boldsymbol{\eta}$ for the dominance case candidates in the $\Delta \delta \boldsymbol{a}$ plane (cf. Fig. \ref{fig: IP RS}, top left) are restricted directionally to the angle of the affine contours, while the $\boldsymbol{\eta}$ in the $\Delta \delta \boldsymbol{e}$ plane (cf. Fig. \ref{fig: IP RS}, bottom left) are defined at the point where the circular $S^*$ contour intersects with $\Delta \delta \boldsymbol{\alpha}$. This means the $\delta e$ dominance case candidate has an infinite number of viable $\boldsymbol{\eta}$, with several possibilities illustrated in Fig. \ref{fig: IP RS}. The IP $\delta v_{\text{min}}$, relevant for hybrid control, are given as
\begin{align} 
    \delta v_{\text{min}, \delta a} &= \pm \frac{1}{2}n_c a_c  \Delta \delta a,  \label{eq: Chernick dv min da}\\
    \delta v_{\text{min}, \delta \lambda} &= \frac{\pm (m \Delta \delta\lambda - \Delta  \delta a) }{ \Delta \delta a_0}, \label{eq: Chernick dv min dl}\\
    \delta v_{\text{min}, \delta e} &= \frac{1}{2} n_c ||a_c \Delta \delta \boldsymbol{e}||_2 ,
    \label{eq: Chernick dv min de}
\end{align}
\noindent where $\Delta \delta a_0$ defines the pseudostate achieved by a unitary tangential impulse at $t_0$, $m = 2/\Phi_{21}^{J_2}(t_0)$, and $\delta v_{\text{min}}$ denotes a propulsive delta-v minimum \cite{chernick}. The authors aligned the $\boldsymbol{\eta}$ from the IP propulsive dominance case candidates with the ROE differential drag reachable sets to find the optimal differential attitude profile associated with each dominance case candidate. The optimal drag maneuvers were identified as maneuvers in both the differential drag $S$ and $S^*$ where the contour of $S^*$ aligns with the hyperplane defined by the propulsive dominance case candidate $\boldsymbol{\eta}$. As a simple example, the single optimal drag maneuver for the $\delta a$ dominance case candidate spans the entire control window in the direction of $\Delta \delta a$ ($\boldsymbol{\eta}_{\delta a}$, cf. Fig. \ref{fig: IP RS}, top right). The analytical effects of the optimal drag profile on each corresponding dominance case candidate's $\delta v_{\text{min}}$ are given as
\begin{align} 
    \Delta \delta v_{\text{min},\delta a} \hspace{-0.5mm} &= \hspace{-0.5mm} \frac{1}{2}a_{c}^{2}n_{c}^{2} \Delta B^{\rho} \Delta t, \label{eq: drag da} \\
    \Delta \delta v_{\text{min},\delta \lambda} \hspace{-0.5mm} &= \hspace{-0.5mm} a_{c}n_c \Delta B^{\rho}\frac{\Hat{\Phi}_{21}^{J_2}m(2 \tau \Delta t \hspace{-0.5mm} - \hspace{-0.5mm} 2t\Delta t \hspace{-0.5mm} - \hspace{-0.5mm} \Delta t^2) \hspace{-0.5mm} - \hspace{-0.5mm} 2\Delta t}{2 \Delta \delta a_0}, \label{eq: drag dl}\\
    \Delta \delta v_{\text{min},\delta e} \hspace{-0.5mm} &= \hspace{-0.5mm} \frac{a_{c}^{2}n_{c}^{2}\Delta B^{\rho}}{\pi} \Delta t, \label{eq: drag de}
\end{align}
\noindent where $\Delta \delta v_{\text{min}}$ denotes the change in the delta-v minimum due to using a drag profile, $\tau = t_f - t_0$, $\Hat{\Phi}_{21}^{J_2}=\Phi_{21}^{J_2}(t_0)/\tau$ is the time invariant component of $\Phi_{21}^{J_2}(t_0)$, $t$ and $\Delta t$ are the start time and length of the dominance case candidate drag profile, and the periodic effects in $\Delta \delta v_{\text{min},\delta e}$ are approximated as linear \cite{hunter_hybrid}. $\Delta \delta v_{\text{min}}$ is negative for all optimal drag maneuvers, although the change may be mathematically defined as positive or negative depending on the sign of $\Delta B^\rho$. These expressions provide a simple analytical evaluation of optimality for reachable set theory-based hybrid control algorithms. They also enable a minimum cost-based estimate on the time required to achieve drag-only controllability of each IP dimension, setting each $\Delta \delta v_{\text{min}}$ to zero and solving for $\Delta t$. This estimate will be used for drag-only reconfiguration maneuver planning within the proposed re-solve architecture.

\section{Re-solve Architecture}

\begin{algorithm} [h]
    \caption{Robust Hybrid Re-Solve Architecture} \label{alg: architecture}
    \begin{algorithmic}[1]
        \Statex $\textbf{Input: } \delta \boldsymbol{\alpha}_f, t_f, \Delta t_{\text{re-solve}}$, DRAG\_ONLY 
        \State $\boldsymbol{\alpha}_c(t_0), \delta \boldsymbol{\alpha}_0, t_{\text{last solve}} = t_0 \gets \text{Navigation} $
        \If{DRAG\_ONLY is true} $t_f \gets t_0,$ (\ref{eq: drag dt}) \EndIf
        \State $\Delta \delta \boldsymbol{\alpha} \gets (\ref{eq: pseudostate}) $
        \State $\Bar{\rho}_c \gets \text{Atmospheric Density Model}$
        \State $t_i, \boldsymbol{u}^{\text{propulsive}}_i, \boldsymbol{u}^{\text{drag}}_i \gets$ Solve (\ref{eq: hybrid control problem}) 
        \While{$t_0 < t_f$}
            \State $\boldsymbol{\alpha}_c(t_0), \delta \boldsymbol{\alpha}_0, t_0 \gets \text{Navigation}$
            \If{$t_0 - t_{\text{last solve}} > \Delta t_{\text{re-solve}}$}
                \State $\Delta \delta \boldsymbol{\alpha}$ $\gets$ (\ref{eq: pseudostate for drag})
                \State Clear $\boldsymbol{u}^{\text{drag}}_i$ for $t_i > t_0$
                \State $\Bar{\rho}_c$ $\gets$ Atmospheric Density Model
                \State $\boldsymbol{u}^{\text{drag}}_i$ $\gets$ Solve (\ref{eq: hybrid control problem})
                \State $t_{\text{last re-solve}} = t_0$
            \EndIf
        \EndWhile
    \end{algorithmic}
\end{algorithm}

\subsection{Hybrid Control}

Algorithm \ref{alg: architecture} (Alg. \ref{alg: architecture}), the hybrid re-solve strategy proposed by this paper, aims to augment the optimality conditions found through reachable set theory with the convergent, error-rejection properties of the iterative methods found in the literature. The architecture begins by using the current absolute and relative state estimates from navigation on-board (Alg. \ref{alg: architecture}, Line 1) to find the reconfiguration pseudostate (Alg. \ref{alg: architecture}, Line 4) and a time-averaged density over the control window from an on-board atmospheric density model (Alg. \ref{alg: architecture}, Line 5) to conduct an initial solve of control at the beginning of the control window (Alg. \ref{alg: architecture}, Line 6) with a generic reachable set theory-based solver. The architecture then periodically re-solves the drag maneuver plan in a diminishing horizon approach, by using the current time and state estimate (Alg. \ref{alg: architecture}, Line 8) to find the remaining pseudostate (Alg. \ref{alg: architecture}, Line 10), removing future drag maneuvers from the maneuver plan while retaining the propulsive maneuvers (Alg. \ref{alg: architecture}, Line 11), updating the time-averaged density (Alg. \ref{alg: architecture}, Line 12), and solving drag maneuvers for the remaining control window (Alg. \ref{alg: architecture}, Line 13), again using a generic reachable set-based solver. The re-solve feedback loop (Alg. \ref{alg: architecture}, Lines 7-16) is restricted to the drag maneuver plan to preserve the cost savings offered by hybrid control, as errors in the dynamics model will likely cause non-optimal maneuvering that translates through propulsive re-solves into increases in overall reconfiguration cost. Drag re-solves are accomplished by removing the modelled effect of future propulsive maneuvers from the current pseudostate (Alg. \ref{alg: architecture}, Line 10), given as
\begin{equation} \label{eq: pseudostate for drag}
    \Delta \delta \boldsymbol{\alpha} = \delta \boldsymbol{\alpha}_f - \boldsymbol{\Phi}^{J_2} (t_0) \delta \boldsymbol{\alpha}_0 - \sum_{i=1}^{k} \boldsymbol{\Gamma}(t_{i}) \boldsymbol{C}^{\text{prop}} \boldsymbol{u}_{i},
\end{equation}
\noindent for $t_i > t_0$ where $\boldsymbol{u}_i$ is the current maneuver plan, produced by either the initial solve in Line 6 or the feedback loop in Line 13, and $t_0$ is set as the current on-board GPS time of the satellite. This approach assumes that the control authority of drag maneuvering is sufficient to reject all modelling errors, a requirement that will be met in later sections. The re-solve cadence used by the architecture is arbitrary, chosen by this paper to be conducted at a constant time interval $\Delta t_{\text{re-solve}}$ over the control window (Alg. \ref{alg: architecture}, Line 9), and may occur as often as allowed by the science mission and the computational constraints of the hardware. The time of the last re-solve is held by $t_{\text{last solve}}$.

\subsection{Differential Drag Only Control}
 
Alg. \ref{alg: architecture} may also accomplish the IP reconfiguration without propellant, similar to the propulsion-free convergent capabilities found in the state of the art. Under drag magnitude limitations, the desired state is no longer guaranteed to be reachable for an arbitrary control window, which must be established before planning optimal maneuvers. A choice of control window that guarantees IP reachability can be made as the sum of control windows required to send all IP dominance case candidate $\delta v_{\text{min}}$ to zero. Setting (\ref{eq: Chernick dv min da}) and (\ref{eq: drag da}) equal, the time $\Delta t_{\delta a}$ required to achieve the $\delta a$ reconfiguration is found as
\begin{align}
    \frac{1}{2}n_c a_c  \Delta \delta a &= \frac{1}{2}a_{c}^{2}n_{c}^{2} \Delta B^{\rho} \Delta t_{\delta a}, \\
    \Delta t_{\delta a} &= \frac{\Delta \delta a}{a_c n_c \Delta B^\rho}.
\end{align}
\noindent The time $\Delta t_{\delta \lambda}$ required to achieve the $\delta \lambda$ reconfiguration is derived by expanding (\ref{eq: Chernick dv min dl}) and (\ref{eq: drag dl}), given as
\begin{align}
    & \delta v_{\text{min},\delta \lambda} = \frac{ \pm (m \Delta \delta\lambda \hspace{-0.5mm} - \hspace{-0.5mm} \Delta  \delta a ) }{ \Delta \delta a_0} \hspace{-0.5mm} \label{eq: Chernick dv min dl expanded}\\
    &= \frac{ \pm \frac{2}{\Hat{\Phi}_{21}^{J_2} \Delta t_{\delta \lambda}} ( \delta \lambda_f - \Hat{\Phi}_{21}^{J_2} \Delta t_{\delta \lambda} \delta a_0 - \delta \lambda_0) - (\delta a_f - \delta a_0) }{ \frac{2}{n_c a_c} }, \nonumber \\
    & \Delta \delta v_{\text{min},\delta \lambda}(t) \hspace{-0.5mm} = \hspace{-0.5mm} a_{c}n_c \Delta B^{\rho}\frac{\Hat{\Phi}_{21}^{J_2}m(2 \tau \Delta t \hspace{-0.5mm} - \hspace{-0.5mm} 2t\Delta t \hspace{-0.5mm} - \hspace{-0.5mm} \Delta t^2) \hspace{-0.5mm} - \hspace{-0.5mm} 2\Delta t}{2 \Delta \delta a_0} \\
    &= \hspace{-0.5mm} a_{c}^{2}n_c^2 \Delta B^{\rho} \Hat{\Phi}_{21}^{J_2}m(\Delta t_{\delta \lambda}^2 \hspace{-0.5mm} - \hspace{-0.5mm} t\Delta t_{\delta \lambda} \hspace{-0.5mm} - \hspace{-0.5mm} 1/4\Delta t_{\delta \lambda}^2) \hspace{-0.5mm} - \hspace{-0.5mm} \Delta t_{\delta \lambda}, \nonumber
\end{align}
\noindent where the length of an optimal $\delta \lambda$ drag maneuver is set to half of the control window $\Delta t_{\delta \lambda}$. Setting (\ref{eq: Chernick dv min dl expanded}) equal to two optimal $\delta \lambda$ maneuvers in the first and second half of the control window, $\Delta t_{\delta \lambda}$ is found as 
\begin{align}
    \delta v_{\text{min},\delta \lambda} = \Delta \delta v_{\text{min},\delta \lambda}&(t = t_0) \\
    &+  \Delta \delta v_{\text{min},\delta \lambda}(t = t_0 + \frac{1}{2} \Delta t_{\delta \lambda}), \nonumber \\
    \frac{1}{4}\Hat{\Phi}_{21}^{J_2} a_c n_c \Delta B^\rho \Delta t_{\delta \lambda}^2 \pm \Hat{\Phi}_{21}^{J_2}&(\delta a_f + \delta a_0)\Delta t_{\delta \lambda} \label{eq: drag dt dl} \\
    &\pm 2(\delta \lambda_0 - \delta \lambda_f) = 0. \nonumber
\end{align}
\noindent $\Delta t_{\delta \lambda}$ is solved as the valid root of (\ref{eq: drag dt dl}). Setting (\ref{eq: Chernick dv min de}) and (\ref{eq: drag de}) equal, the time $\Delta t_{\delta e}$ guaranteed to achieve the $\delta e$ reconfiguration is given as
\begin{align}
    \frac{1}{2} n_c ||a_c \Delta \delta \boldsymbol{e}||_2 &= \frac{a_{c}^{2}n_{c}^{2}\Delta B^{\rho}}{\pi} \Delta t_{\delta e}, \\
    \Delta t_{\delta e} &= \frac{\pi ||\Delta \delta \boldsymbol{e}||_{2,\text{max}}}{2a_c n_c \Delta B^\rho},
\end{align}
\noindent where $||\Delta \delta \boldsymbol{e}||_{2,\text{max}} = ||\delta \boldsymbol{e}_f||_2+||\delta \boldsymbol{e}_0||_2$ is the largest possible relative eccentricity pseudostate. $||\Delta \delta \boldsymbol{e}||_{2,\text{max}}$ is used to find a closed-form expression for $||\Delta \delta \boldsymbol{e}||_2$, which is otherwise analytically intractable due to time-dependent trigonometric functions in $\boldsymbol{\Phi}^{J_2}(t)$. Therefore, a length of control window $t_f-t_0$ that guarantees IP reachability can be found a-priori as
\begin{equation} \label{eq: drag dt}
    t_f-t_0 = 2(\Delta t_{\delta a} + \Delta t_{\delta \lambda} + \Delta t_{\delta e}),
\end{equation}
\noindent The sum $\Delta t_{\delta a} + \Delta t_{\delta \lambda} + \Delta t_{\delta e}$ is doubled to prevent navigation errors from producing a control window that is too short for IP drag-only reachability. This "2" coefficient is arbitrary and can be tuned lower or higher depending on the expected performance of state estimation. For drag-only control, the proposed architecture modifies the length of the control window to ensure the desired state is reachable using drag-only maneuver planning (Alg. \ref{alg: architecture}, Line 2), and the drag-only mode is initiated by setting the boolean flag DRAG\_ONLY to true.

\section{Provably Lyapunov-Stable Control}

\subsection{Ideal Control}

Online control methods commonly define Lyapunov stability in terms of tracking error for the single, current maneuver to converge to a desired trajectory. The hybrid control formulation in (\ref{eq: hybrid control problem}) allows a much simpler interpretation of stability, converging directly to the desired state or pseudostate-space origin. However, complete maneuver planning requires stability guarantees for all future planned maneuvers, instead of just the current maneuver. This section will establish the Lyapunov stability of reachable set theory-based hybrid control (Alg. \ref{alg: architecture}, Lines 6 and 13) under ideal and optimal actuation and dynamics. The decoupling of drag and propulsive control resulting from the reachable set theory approach allows each component of the hybrid solver to be analyzed separately.

Consider a Lyapunov function $V$ proportional to the square of the IP dominance case $\delta v_{\text{min}}$ numerator in (\ref{eq: delta-v min}), defined as
\begin{equation} \label{eq: V}
    V = \frac{1}{2} \left ( \boldsymbol{\eta}^{\text{T}} \Delta \delta \boldsymbol{\alpha} \right )^{\text{T}} \left ( \boldsymbol{\eta}^{\text{T}} \Delta \delta \boldsymbol{\alpha} \right ),
\end{equation}
\noindent where $\boldsymbol{\eta}$ is the outward normal of the reconfiguration dominance case intersecting $\Delta \delta \boldsymbol{\alpha}$ at optimal cost as illustrated in Fig. \ref{fig: IP RS}. $V$ is radially unbounded, positive definite, continuously differentiable, and, as a function of the reconfiguration pseudostate, remains constant when the deputy is not maneuvering. Additionally, $\boldsymbol{\eta}$ is property of the reconfiguration and, like the pseudostate, is invariant. This implies that no optimal maneuver will cause the reconfiguration dominance case to change to another dominance case candidate, such that the $\delta v_{\text{min}}$ for the reconfiguration dominance case will always be greater than the $\delta v_{\text{min}}$ for the other dominance case candidates throughout the reconfiguration. Therefore, decreasing the reconfiguration dominance case $\delta v_{\text{min}}$ to zero must also decrease the other dominance case candidate $\delta v_{\text{min}}$ to zero, and the entire reconfiguration is completed even if $\boldsymbol{\eta}$ contains one or more zeros. The change in $V$ due to an optimal maneuver $\boldsymbol{u}_i$ at time $t_i$ is given from the dynamic constraint in (\ref{eq: hybrid control problem}) as
\begin{align} \label{eq: V_dot}
    \Dot{V} &= \left ( \boldsymbol{\eta}^{\text{T}} \Delta \delta \boldsymbol{\alpha} \right )^{\text{T}} \left ( \boldsymbol{\eta}^{\text{T}} \Dot{\Delta \delta \boldsymbol{\alpha}} \right ) \\
    &= \left ( \boldsymbol{\eta}^{\text{T}} \Delta \delta \boldsymbol{\alpha} \right )^{\text{T}} \left ( \boldsymbol{\eta}^{\text{T}} (\Delta \delta \boldsymbol{\alpha} - \boldsymbol{\Gamma}(t_i)\boldsymbol{u}_i) - \boldsymbol{\eta}^{\text{T}} \Delta \delta \boldsymbol{\alpha}  \right ) \nonumber \\
    &= - \left ( \boldsymbol{\eta}^{\text{T}} \Delta \delta \boldsymbol{\alpha} \right )^{\text{T}} \left ( \boldsymbol{\eta}^{\text{T}} \boldsymbol{\Gamma}(t_i)\boldsymbol{u}_i \right ). \nonumber
\end{align}
\noindent Koenig and D'Amico showed that $\boldsymbol{\eta}^{\text{T}} \Delta \delta \boldsymbol{\alpha} > 0$ for an optimal $\boldsymbol{\eta}$, as the minimum cost must be positive \cite{koenig_RST}. For optimal propulsion, maneuvers will always act in the pseudostate direction defined by $\boldsymbol{\eta}$ to decrease and minimize the remaining IP reconfiguration cost. Therefore, $\boldsymbol{\eta}^{\text{T}} \boldsymbol{\Gamma}(t_i)\boldsymbol{u}_i > 0$ and $\Dot{V}$ is negative definite. The origin of pseudostate space, and the desired state at the end of the control window, is globally asymptotically Lyapunov stable for ideal optimal propulsive control. 

For optimal differential drag, the offline maneuver planning approach of the reachable set solver (Alg. \ref{alg: architecture}, Lines 6 and 13) ensures that the differential drag maneuver plan in its entirety optimally minimizes the IP propulsive minimum cost of the reconfiguration. However, unlike propulsive maneuvers, each individual differential drag maneuver is not required to greedily decrease the remaining IP reconfiguration cost as a condition of optimality, as drag is a zero or negligible cost maneuver and is not included in the cost function in (\ref{eq: hybrid control problem}). Therefore, $\Dot{V}$ is negative semi-definite for ideal optimal differential drag control. 

Asymptotic stability can be verified instead through LaSalle's Invariance Principle Theorem \cite{lasalle}. Let $D$ be the set of pseudostates within a $\gamma$-level set of $V$, defined as
\begin{equation}
    D = \{\Delta \delta \boldsymbol{\alpha}: V(\Delta \delta \boldsymbol{\alpha}) \leq \gamma \}.
\end{equation}
\noindent It is worth noting that the reachable sets defined by $S^*$ in (\ref{eq: RST S*}) are also level sets of the linearly scaling cost $c$, and $D$ is almost equivalent to $S^*$, having the same shape but defined instead by the proportional square of minimum cost. $\Dot{V}$ from (\ref{eq: V_dot}) is negative semi-definite under ideal optimal drag maneuvering for all $\Delta \delta \boldsymbol{\alpha} \in D$, such that $D$ is forward invariant. Let $E$ be the set of all states in $D$ where $\Dot{V} = 0$, given as
\begin{align}
    E = \{ &\Delta \delta \boldsymbol{\alpha} \in D, \boldsymbol{u}_i: \Dot{V} = 0 \} \\
    = \{ &\Delta \delta \boldsymbol{\alpha} = 0, \boldsymbol{u}_i = 0, \boldsymbol{\eta}^{\text{T}} \boldsymbol{\Gamma}(t_i)\boldsymbol{u}_i = 0 \}. \nonumber
\end{align}
\noindent $E$ contains the origin when the reconfiguration is achieved, the non-maneuvering deputy, and maneuvers that do not change the IP dominance case $\delta v_{\text{min}}$. Let $N$ be the largest set of invariant states in $E$, given as
\begin{equation}
    \hspace{-1mm} N = \{ \Delta \delta \boldsymbol{\alpha} \in E: \Dot{\Delta \delta \boldsymbol{\alpha}} = 0 \} = \{ \Delta \delta \boldsymbol{\alpha} = 0, \boldsymbol{u}_i = 0 \}. \hspace{-1mm}
\end{equation}
\noindent $N$ contains the achieved reconfiguration and the non-maneuvering spacecraft. Therefore, the origin of pseudostate space is globally asymptotically stable under ideal optimal differential drag control.

\subsection{Error Bounds on Lyapunov Stability}

The purpose of Alg. \ref{alg: architecture} is to reject all forms of real-world error that deteriorate hybrid control performance. The stability requirements within Lyapunov's Second Method and LaSalle's Invariance Principle Theorem enable convergence bounds to be placed on the error relative to the optimal maneuver plan and provide intuition on the best dynamic modeling feedback approach (Alg. \ref{alg: architecture}, Lines 7-16). Consider an arbitrary general positive scalar error term $\epsilon$ for any non-ideal maneuver $\Tilde{\boldsymbol{u}}_i$  in the hybrid maneuver plan, such that $\boldsymbol{\eta}^{\text{T}} \boldsymbol{\Gamma}(t_i)  \Tilde{\boldsymbol{u}}_i = \boldsymbol{\eta}^{\text{T}} \boldsymbol{\Gamma}(t_i)  \boldsymbol{u}_i \pm \epsilon$. Throughout this paper, quantities denoted with a tilde are associated with realized in-flight values under typical operational errors corresponding to ideal non-tilde values used on-board in maneuver planning. In this section, $\epsilon$ captures both maneuver and dynamic modeling error in a single term. Substituting the non-ideal maneuver into $\Dot{V}$ in (\ref{eq: V_dot}) results in
\begin{equation} \label{eq: V_dot err}
    \Dot{V} = \left ( \boldsymbol{\eta}^{\text{T}} \Delta \delta \boldsymbol{\alpha} \right )^{\text{T}} \left( \pm \epsilon - \boldsymbol{\eta}^{\text{T}} \boldsymbol{\Gamma}(t_i)  \Tilde{\boldsymbol{u}_i} \right).
\end{equation}
\noindent The error bound for stable convergence can be established by examining the effect of $\epsilon$ on the change in optimal cost, defined as
\begin{align}
    \boldsymbol{\eta}^{\text{T}} \Dot{\Delta \delta \boldsymbol{\alpha}} = &\pm \epsilon -\boldsymbol{\eta}^{\text{T}} \boldsymbol{\Gamma}(t_i)  \Tilde{\boldsymbol{u}}_i \leq 0, \\
    \pm \epsilon &\leq \boldsymbol{\eta}^{\text{T}} \boldsymbol{\Gamma}(t_i)  \Tilde{\boldsymbol{u}}_i .
\end{align}
\noindent Two error bounds can be derived from the positive and negative errors by resubstituting the definition $\pm \epsilon = \boldsymbol{\eta}^{\text{T}} \boldsymbol{\Gamma}(t_i)  \Tilde{\boldsymbol{u}}_i - \boldsymbol{\eta}^{\text{T}} \boldsymbol{\Gamma}(t_i)  \boldsymbol{u}_i$. $+\epsilon$ provides a lower bound for convergence, given as
\begin{align} 
    \epsilon = \boldsymbol{\eta}^{\text{T}} \boldsymbol{\Gamma}(t_i)  \Tilde{\boldsymbol{u}}_i &- \boldsymbol{\eta}^{\text{T}} \boldsymbol{\Gamma}(t_i)  \boldsymbol{u}_i \leq \boldsymbol{\eta}^{\text{T}} \boldsymbol{\Gamma}(t_i)  \Tilde{\boldsymbol{u}}_i, \\
    &\boldsymbol{\eta}^{\text{T}} \boldsymbol{\Gamma}(t_i)  \boldsymbol{u}_i \geq 0, \label{eq: V_dot err low}
\end{align}
\noindent while $-\boldsymbol{\epsilon}$ provides an upper bound for convergence, given as
\begin{align} 
    -\epsilon = -(\boldsymbol{\eta}^{\text{T}} \boldsymbol{\Gamma}&(t_i)  \Tilde{\boldsymbol{u}}_i - \boldsymbol{\eta}^{\text{T}} \boldsymbol{\Gamma}(t_i)  \boldsymbol{u}_i) \leq \boldsymbol{\eta}^{\text{T}} \boldsymbol{\Gamma}(t_i)  \Tilde{\boldsymbol{u}}_i, \\
    &\boldsymbol{\eta}^{\text{T}} \boldsymbol{\Gamma}(t_i)  \boldsymbol{u}_i \leq 2\boldsymbol{\eta}^{\text{T}} \boldsymbol{\Gamma}(t_i)  \Tilde{\boldsymbol{u}}_i, \label{eq: V_dot err up}
\end{align}
\noindent and combining (\ref{eq: V_dot err low}) and (\ref{eq: V_dot err up}) results in the constraint
\begin{equation} \label{eq: converge bounds}
    0 \leq \boldsymbol{\eta}^{\text{T}} \boldsymbol{\Gamma}(t_i)  \boldsymbol{u}_i \leq 2\boldsymbol{\eta}^{\text{T}} \boldsymbol{\Gamma}(t_i)  \Tilde{\boldsymbol{u}}_i.
\end{equation}
\noindent In words, the proposed re-solve architecture, Alg. \ref{alg: architecture}, will asymptotically converge if the on-board modelled optimal maneuver effect on the IP minimum cost in the generic reachable set solver (Alg. \ref{alg: architecture}, Lines 6 and 13) is no more than double the effect of the actual, realized maneuver. Relating to traditional control techniques, the lower bound in (\ref{eq: V_dot err low}) is analogous to underdamped control, with iterative overshoot corrections over the control window, and the upper bound in (\ref{eq: V_dot err up})  is analogous to a bound on overdamping, where the re-solve architecture will no longer approach the desired final state due to undershoot.

\section{Primary Sources of Modeling Error}

The re-solve architecture is designed to meet the stability constraint posed in (\ref{eq: converge bounds}) by modifying the inputs or resulting maneuver plan outputs of the hybrid control problem (Alg. \ref{alg: architecture}, Lines 6 and 13), instead of requiring a specialized reachable set-based solver with a unique derivation of Lyapunov stability that loses generality. Dynamic modeling inaccuracies are produced in-flight from a variety of solver input and output sources, including uncertain navigation, non-ideal hardware, highly variable environmental parameters, and the LTV model of relative spacecraft motion itself. The primary sources of error for hybrid control are listed here and organized into two categories, linear and nonlinear, based on their effects on the output final ROE state at the end of the control window. These categories will determine how the various errors are analyzed in further sections. 

The propulsive linear errors are magnitude errors, driven by thruster imperfections and differences between ground testing and flight performance. Linear differential drag errors are also magnitude errors caused by the augmented differential ballistic coefficient $\Delta B^\rho$. Errors in $\Delta B^\rho$ are primarily produced from inaccurate atmospheric density modelling, but also can result from poor characterization of the spacecraft's coefficient of drag or poor attitude determination and control of the spacecraft's cross-sectional area. Both types of magnitude error cause a linear scaling of the resulting reachable set for the maneuver, as a change of the linear variable $\boldsymbol{u}_i$ in the constraint of (\ref{eq: hybrid control problem}). Navigation uncertainty in the initial relative state between the chief and deputy effects the pseudostate of the reconfiguration problem, as a change of the linear variable $\delta \boldsymbol{\alpha}_0$ in (\ref{eq: pseudostate}), and causes a linear translation of the desired state in pseudostate-space.

For nonlinear error effects, inaccurate maneuver timing affects both propulsive and differential drag control through the translation of software commands to hardware actuation. Propulsive timing delays can also result from the ramp-up or activation of the thruster as well as differences between the impulsive-approximated maneuver and the actual thrust curve achieved by the thruster. Differential drag timing delays can be produced by the difference between the instantaneous-approximated attitude changes and the actual rotational dynamics of the spacecraft. 
Directional inaccuracies also cause nonlinear effects, with propulsive errors caused by attitude determination and control and differential drag errors resulting from the shared RTN frame approximation between the chief and deputy. All of these nonlinear maneuvering effects, excluding errors in propulsive maneuver direction, cause the resulting reachable set to warp in shape. Propulsive direction errors, holding maneuver cost and time constant, will cause no change in reachable set shape, as all variations of viable maneuver direction are already captured in the creation of the reachable sets. Navigation uncertainty in the absolute Keplerian state nonlinearly affects both the reconfiguration pseudostate and the effect of each maneuver, causing both a translation of the desired state in pseudostate-space and a warping of the reachable sets. 

Finally, two inherent flaws exist in the previous Lyapunov analysis. Closed-form maneuver planning algorithms based on reachable set theory are not always guaranteed to output an optimal maneuver plan, given potential dynamic coupling between the 2D planes (i.e. IP/OOP ROE coupling). However, the solver will find the lowest cost within the closed-form approach to ensure near-optimality for all output maneuver plans. These plans can be considered "approximations" of the true optimal maneuver plan, caused by magnitude errors that linearly scale the reachable sets and nonlinear differences in maneuver direction and timing within the solver. Additionally, the Lyapunov formulation uses the reconfiguration pseudostate as its state vector, implying the ROE LTV dynamics model perfectly captures relative spacecraft motion. This is not the case, as the full model is loses accuracy at larger separations and neglects lower order orbital perturbations to Keplerian motion, producing small, minor inaccuracies. These approximations allow generic reachable set-based solvers of (\ref{eq: hybrid control problem}) within Alg. \ref{alg: architecture} to be covered by the same Lyapunov formulation. Both resulting errors are difficult to quantify and decouple from other parameters of the control problem, such as initial and final state, the length of the control window, and the exact reachable set-based solver. However, neither error source is expected to significantly effect the performance of the re-solve architecture, and they are not considered in the following error covariance analysis.

\section{Covariance Sensitivity Analysis}

The uncertainty of each error source can be quantified and translated to the achieved final state using a covariance sensitivity analysis to produce experimental performance metrics. Comparing these performance metrics will inform the robust control strategy of Alg. \ref{alg: architecture}, by conditioning the inputs and outputs of the reachable set-based solver (Alg. \ref{alg: architecture}, Lines 6 and 13) with buffer functions (i.e. scale factor, min/max bounds) that ensure stable convergence. This section extends the analysis approach proposed by Chernick and D'Amico for propulsive control to hybrid control \cite{chernick}. Matching the distinction made while identifying the primary sources of error in hybrid control, linear error effects are derived analytically, while quantifying nonlinear error effects requires numerical methods. Throughout this section, arbitrary errors are given as a scalar $\epsilon$ or vector $\boldsymbol{\epsilon}$ to produce the non-ideal maneuver planning quantity, denoted as in previous sections by a tilde. The source of the error is identified by the subscript of $\epsilon$. All error sources are modelled as a zero mean multi-variate Gaussian distribution with diagonal covariance matrix $\boldsymbol{COV}_\epsilon = diag(\sigma_1,\sigma_2,...,\sigma_n)$ for an $n$D state with associated standard deviations $\sigma_i$.

Analytical expressions for the linear effects are derived by separating the resulting final state into ideal and error components. Considering a general maneuver magnitude error formulation $\Tilde{\boldsymbol{u}}_{||.||_2} = \boldsymbol{u} + \boldsymbol{\epsilon}_{|\boldsymbol{u}|} = \boldsymbol{u} + \epsilon_{|\boldsymbol{u}|} \Hat{\boldsymbol{u}}$, where $\Hat{\boldsymbol{u}}$ is the normalized direction of $\boldsymbol{u}$, and initial relative state error formulation $\Tilde{\delta \boldsymbol{\alpha}}_0 = \delta \boldsymbol{\alpha}_0 + \boldsymbol{\epsilon}_{\delta \boldsymbol{\alpha}_0}$, the dynamics model from (\ref{eq: hybrid control problem}) is expanded and formulated as
\begin{equation} \label{eq: linear err dynamics}
    \delta \boldsymbol{\alpha}_f = \boldsymbol{\Phi}^{J_2} (t_0) \Tilde{\delta \boldsymbol{\alpha}_0} + \sum_{i=1}^{k} \boldsymbol{\Gamma}(t_i) \Tilde{\boldsymbol{u}}_{||.||_2,i}.
\end{equation}
\noindent This expression is split into the ideal final state $\delta \boldsymbol{\alpha}_f^{\text{ideal}}$ and final state error $\delta \boldsymbol{\alpha}_f^\epsilon$, given as
\begin{align}
    &\delta \boldsymbol{\alpha}_f = \delta \boldsymbol{\alpha}_f^{\text{ideal}} + \delta \boldsymbol{\alpha}_f^\epsilon, \\
    &\delta \boldsymbol{\alpha}_f^{\text{ideal}} = \boldsymbol{A}^{\delta \boldsymbol{\alpha}_0} \delta \boldsymbol{\alpha}_0 + \boldsymbol{A}^{\text{prop}}\boldsymbol{u}_i + \boldsymbol{A}^{\text{drag}}\boldsymbol{u}_i, \\
    &\delta \boldsymbol{\alpha}_f^\epsilon = \boldsymbol{A}^{\delta \boldsymbol{\alpha}_0} \boldsymbol{\epsilon}_{\delta \boldsymbol{\alpha}_0} + \boldsymbol{A}^{\text{prop}}\boldsymbol{\epsilon}_{||\boldsymbol{u}^{\text{prop}}_i||_2} + \boldsymbol{A}^{\text{drag}}\boldsymbol{\epsilon}_{||\boldsymbol{u}^{\text{drag}}_i||_2},
\end{align}
\noindent with the substitutions $\boldsymbol{A}^{\delta \boldsymbol{\alpha}_0} = \boldsymbol{\Phi}^{J_2} (t_0)$, $\boldsymbol{A}^{\text{prop}} = \sum_{i=1}^{k} \boldsymbol{\Gamma}(t_i) \boldsymbol{C}^{\text{prop}}$ and $\boldsymbol{A}^{\text{drag}} =  \sum_{i=1}^{k} \boldsymbol{\Gamma}(t_i) \boldsymbol{C}^{\text{drag}}$ where $\boldsymbol{C}^{\text{drag}} = [\boldsymbol{0}_3,\boldsymbol{I}_3]$. $\delta \boldsymbol{\alpha}_f$ is a weighted sum of different Gaussian distributions and, assuming that the errors considered by this analysis are independent, is itself a Gaussian with mean $\delta \boldsymbol{\alpha}_f^{\text{ideal}}$. The resulting uncorrelated effect of each error source on the $\delta \boldsymbol{\alpha}_f$ can be found independently by linearly transforming the covariance of the error source, given as
\begin{equation} \label{eq: linear error final state cov}
    \boldsymbol{COV}_{\delta \boldsymbol{\alpha}_f} = \boldsymbol{A}[\boldsymbol{COV}_\epsilon]\boldsymbol{A}^{\text{T}},
\end{equation}
\noindent with the applicable $\boldsymbol{A}$ for each error source.

Analytical expressions are not possible for the nonlinear covariance effects, requiring the use of numerical methods. Considering a general maneuver timing error formulation $\Tilde{t} = t + \epsilon_{t}$, maneuver direction error formulation $\Tilde{\boldsymbol{u}}_\measuredangle = \boldsymbol{R}_z(\epsilon_{\measuredangle \boldsymbol{u},1})\boldsymbol{R}_x(\epsilon_{\measuredangle \boldsymbol{u},2})\boldsymbol{R}_z(\epsilon_{\measuredangle \boldsymbol{u},3})\boldsymbol{u}$ using a 3-1-3 Euler angle rotation with rotation matrix $\boldsymbol{R}$ and three distinct error angles, and initial absolute state error formulation $\Tilde{\boldsymbol{\alpha}}_c = \boldsymbol{\alpha}_c + \boldsymbol{\epsilon}_{\boldsymbol{\alpha}_c}$, the resulting dynamics model is formulated as
\begin{equation} \label{eq: nonlinear err dynamics}
    \delta \boldsymbol{\alpha}_f = \Tilde{\boldsymbol{\Phi}}^{J_2} (t_0) \delta \boldsymbol{\alpha}_0 + \sum_{i=1}^{k} \Tilde{\boldsymbol{\Gamma}}(\Tilde{t_{i}}) \Tilde{\boldsymbol{u}}_{\measuredangle,i},
\end{equation}
\noindent where $\Tilde{\boldsymbol{\Gamma}}(t) = [\Tilde{\boldsymbol{\Phi}}^{J_2}(t) \Tilde{\boldsymbol{B}}(t),\int_{t_{i}}^{t_{i+1}}\Tilde{\boldsymbol{\Phi}}^{J_2}(t)\Tilde{\boldsymbol{B}}(t)dt]$, $\Tilde{\boldsymbol{\Phi}}^{J_2} (t) = \boldsymbol{\Phi}^{J_2} (\Tilde{\boldsymbol{\alpha}}_c (t), t_f-t)$, and $\Tilde{\boldsymbol{B}}(t) = \boldsymbol{B}(\Tilde{\boldsymbol{\alpha}}_c (t))$. Again assuming that the errors considered by this analysis are independent, Gaussian Monte Carlo samples are taken for a single nonlinear error source and converted to $\delta \boldsymbol{\alpha}_f$. By law of large numbers, the covariance of a sufficiently large number of samples will converge to the true $\boldsymbol{COV}_{\delta \boldsymbol{\alpha}_f}$ for that error source. 

The on-diagonal terms of $\boldsymbol{COV}_{\delta \boldsymbol{\alpha}_f}$ are the independent squared standard deviation of each ROE, but in the presence of non-zero off-diagonal terms, standard deviation is insufficient to establish performance metrics or confidence bounds. For completeness, the uncertainty coupling effects between the ROE need to be quantified. The Gaussian error model is leveraged to find a $n$D confidence bound ellipsoid, defined in shape by the covariance matrix $\boldsymbol{COV}_{\delta \boldsymbol{\alpha}_f}$ and  scaled by the desired critical value of the chi-square cumulative distribution function. The performance bounds are taken at the $n$D ellipsoid maxima, the maximum projection of the ellipsoid on each dimension, as a worst-case conservative value. As stated by Chernick and D'Amico, this can be interpreted as fitting an $n$D polygon, representing the performance metrics, to the desired confidence bound $n$D ellipsoid such that the polygon is centered at the mean, defined here as $\delta \boldsymbol{\alpha}_f^{\text{ideal}}$, and tangent to the $n$D ellipsoid at its maxima \cite{chernick}.

\begin{figure}
    \centering
    \includegraphics[scale=0.23]{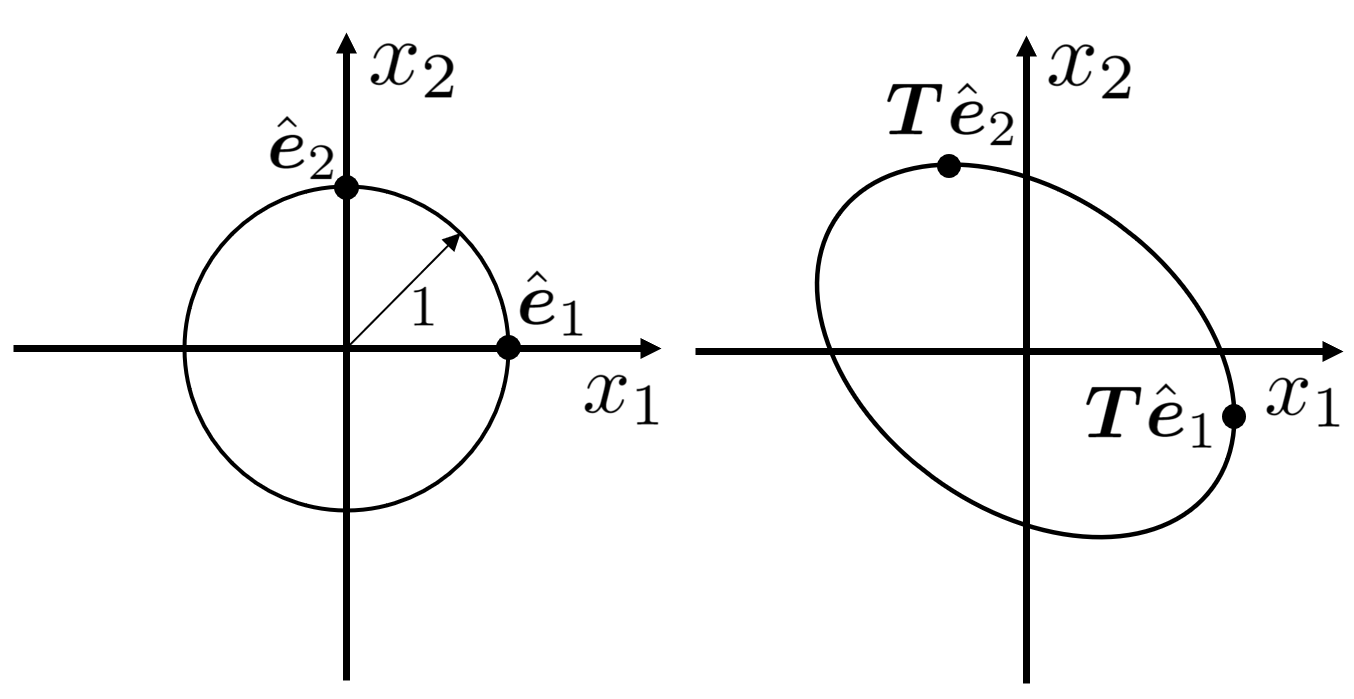}
    \caption{Points on a 2-dimensional unit sphere centered on the origin (left) can be translated to a 2-dimensional ellipsoid centered on the origin (right) using transformation matrix $T$. This transformation is shown for points on the unit sphere that coincide with the $x_1$ and $x_2$ axes.}
    \label{fig: sphere to ellipsoid}
\end{figure}

To find these ellipsoid maxima, first consider an $n$D unit spheroid. Any point on this unit spheroid can be linearly transformed to an equivalent point on a $n$D ellipsoid, such that the transformation matrix $\boldsymbol{T}$ between a covariance ellipsoid and unit spheroid is defined as $\boldsymbol{COV}_{\delta \boldsymbol{\alpha}_f} = \boldsymbol{T}\boldsymbol{T}^{\text{T}}$ \cite{chernick}. Fig. \ref{fig: sphere to ellipsoid} illustrates this transformation for a 2D plane, using coordinates $x_1$ and $x_2$ and transforming the points on the unit sphere that match the locations of the coordinate basis vectors $\Hat{\boldsymbol{e}}_1$ and $\Hat{\boldsymbol{e}}_2$.

\begin{figure}
    \centering
    \includegraphics[scale=0.28]{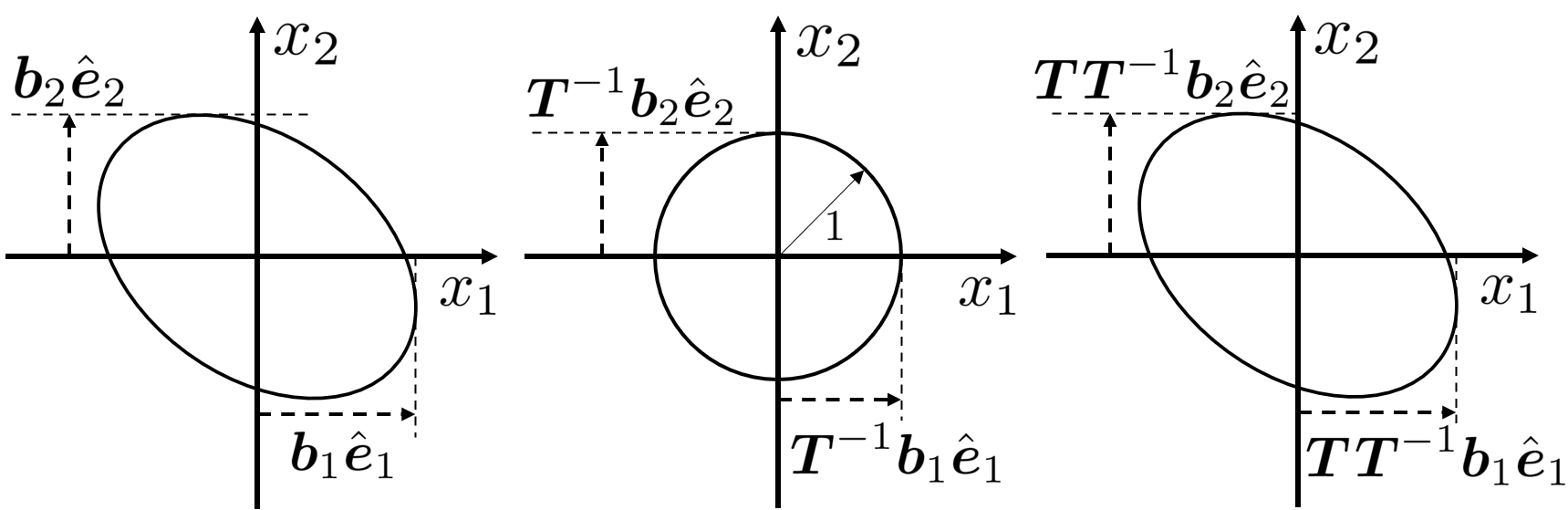}
    \caption{Performance metrics $b_i = \boldsymbol{b}_i\hat{\boldsymbol{e}}_i$ on 2-dimensional confidence bound ellipsoid (left), inversely transformed to points on a unit sphere (center), and transformed back to the confidence bound ellipsoid (right).}
    \label{fig: perform bound formulation}
\end{figure}

Given that a valid covariance matrix is positive semi-definite and full rank, the eigendecomposition of the covariance matrix results in diagonal eigenvalue matrix $\boldsymbol{\lambda}$ and eigenvector matrix $\boldsymbol{v}$ where column $i$ in $\boldsymbol{v}$ is the eigenvector $\boldsymbol{v}_i$ corresponding to eigenvalue $\lambda_i$ in column $i$ and row $i$ of $\boldsymbol{\lambda}$. The eigendecomposition can be rearranged to realize the transformation matrix from the unit spheroid to the covariance ellipsoid, derived as
\begin{equation}
    \boldsymbol{COV}_{\delta \boldsymbol{\alpha}_f} = \boldsymbol{v}\boldsymbol{\lambda}\boldsymbol{v}^{\text{T}} = \boldsymbol{v}\boldsymbol{\lambda}^{1/2}[\boldsymbol{v}\boldsymbol{\lambda}^{1/2}]^{\text{T}} = \boldsymbol{T}\boldsymbol{T}^{\text{T}}.
\end{equation}
\noindent Now consider scaling $\boldsymbol{COV}_{\delta \boldsymbol{\alpha}_f}$ and the corresponding $\boldsymbol{T}$ by the chi-square critical value $\chi^2$ for the desired cumulative probability and the degrees of freedom of the state, resulting in an ellipsoid that represents the desired confidence bound. The transformation matrix from a unit spheroid to the confidence bound ellipsoid is given as
\begin{equation}
    \boldsymbol{T} = 2\chi\boldsymbol{v}\boldsymbol{\lambda}^{1/2}.
\end{equation}
\noindent As worst-case conservative values, each performance metric $b_i = \boldsymbol{b}_i\hat{\boldsymbol{e}}_i$ corresponds to a hyperplane defined by the basis vector $\hat{\boldsymbol{e}}_i$ of the corresponding dimension and tangent to the confidence bound ellipsoid at its maxima $\boldsymbol{b}_i$. Fig. \ref{fig: perform bound formulation} visually illustrates the following process used to formulate these $b_i$ analytically.

Focusing on the alignment with the basis vector, consider transforming an ellipsoid maxima to its corresponding point on the unit spheroid $\boldsymbol{T}^{-1}\boldsymbol{b}_i\hat{\boldsymbol{e}}_i$. This point can equivalently be expressed as the normalized column $i$ of $\boldsymbol{T}^{-1}$ because this column completely defines the direction of transformed maximum on the unit spheroid. Therefore, the performance bound on each dimension $i$ can be easily found by transforming this simplified expression of the spheroid point back to a point on the ellipsoid, given as
\begin{equation} \label{eq: performance bound}
    b_i = \frac{\boldsymbol{T}\boldsymbol{T}^{-1}\hat{\boldsymbol{e}}_i}{||\boldsymbol{T}^{-1}\hat{\boldsymbol{e}}_i||_2}.
\end{equation}

\section{Hybrid Error Analysis}

The covariance uncertainty analysis detailed in the previous section is conducted for three different reconfigurations. Each reconfiguration involves a different IP dominance case candidate to ensure the analysis covers all variations of the optimal differential drag profile. These three reconfigurations are defined for a 100 orbit control window, given as
\begin{align}
    \boldsymbol{\alpha}_c (t_0) &= 
    \begin{bmatrix}
        a & e & i & \Omega & \omega & M
    \end{bmatrix} \\
    &= 
    \begin{bmatrix}
        6798\text{km} & 0.001 & 51^\circ & 30^\circ & 200^\circ & 45^\circ
    \end{bmatrix}, \nonumber \\
    \begin{bmatrix}
        a_c \delta \boldsymbol{\alpha}_f \\ a_c \delta \boldsymbol{\alpha}_{0,\delta a} \\
        a_c \delta \boldsymbol{\alpha}_{0,\delta \lambda} \\ a_c \delta \boldsymbol{\alpha}_{0,\delta e}
    \end{bmatrix} \hspace{-1mm} &= \hspace{-1mm}
    \begin{bmatrix}
        0 & 5 & 0 & 1 & 0 & 1 \\
        0.8 & \hspace{-1mm} 150 & \hspace{-1mm} 0.55 & \hspace{-1mm} 0.85 & \hspace{-1mm} -0.05 & \hspace{-1mm} 0.95 \\
        -0.02 & \hspace{-1mm} 100 & \hspace{-1mm} 0.5 & \hspace{-1mm} 0.9 & \hspace{-1mm} 0.05 & \hspace{-1mm} 0.95 \\
        0.02 & \hspace{-1mm} 50 & \hspace{-1mm} 0.1 & \hspace{-1mm} 1 & \hspace{-1mm} -0.05 & \hspace{-1mm} 0.95 \\
    \end{bmatrix} \hspace{-0.5mm} \text{km} ,
\end{align}
\noindent where subscripts $\delta a, \delta \lambda,$ and $\delta e$ denote initial conditions for the associated dominance case. The generic reachable set-based solver (Alg. \ref{alg: architecture}, Lines 6 and 13) chosen to produce maneuver plans for this analysis is a hybrid solver designed in previous work by the authors, which produces provably-optimal hybrid control solutions in closed-form \cite{hunter_hybrid}. To give context on the balance between propulsive and drag control in these reconfigurations, the $\delta a$ dominant reconfiguration saves 0.1323m/s of $\delta v$ cost with hybrid control over propulsion alone, solved from (\ref{eq: drag da}), and requires a total IP $\delta v$ of 0.4506m/s, found by subtracting (\ref{eq: drag da}) from (\ref{eq: Chernick dv min da}). The $\delta \lambda$ dominant reconfiguration saves 0.0666m/s of $\delta v$ cost with hybrid control over propulsion alone, solved from (\ref{eq: drag dl}), and requires a total IP $\delta v$ of 0.0592m/s, found by subtracting (\ref{eq: drag dl}) from (\ref{eq: Chernick dv min dl}). The $\delta e$ dominant reconfiguration saves 0.0831m/s of $\delta v$ cost with hybrid control over propulsion alone, solved from (\ref{eq: drag de}), and requires a total IP $\delta v$ of 0.1083m/s, found by subtracting (\ref{eq: drag de}) from (\ref{eq: Chernick dv min de}). As mentioned previously, each error source will be modelled as a zero-mean Gaussian, and the standard deviations $\sigma$ for each error source are given in Table \ref{tab: error sources}.
\begin{table}[h]
    \centering
    \setlength{\tabcolsep}{4pt}
    \caption{Standard deviations associated with the Gaussian model of each hybrid error source. Superscripts "prop" and "drag" are added to the associated maneuver errors.}
    \begin{tabular}{c|c||c|c}
        \multicolumn{2}{c}{Linear Sources} & \multicolumn{2}{c}{Nonlinear Sources} \\
        \hline \hline
        $\epsilon$ & $\sigma^2$ & $\epsilon$ & $\sigma^2$ \\
        \hline
        $|\boldsymbol{u}^{\text{prop}}_i|$ & $0.01^2$ (m/s)$^2$ & $t_i^{\text{prop}}, t_j^{\text{drag}}$ & $60^2$ s$^2$ \\
        $|\boldsymbol{u}^{\text{drag}}_j|$ & $(5 \times 10^{-12})^2$ (m/s)$^2$ & $\measuredangle  \boldsymbol{u}^{\text{prop}}_i$ & $1^2$ ($^\circ$)$^2$\\
         & & $\measuredangle  \boldsymbol{u}^{\text{drag}}_j$ & $5^2$ ($^\circ$)$^2$ \\
        $\delta x_0$ & $0.01^2$ m$^2$ & $x_0$ & $1^2$ m$^2$ \\
        $\delta \Dot{x}_0$ & $0.001^2$ (m/s)$^2$ & $\Dot{x}_0$ & $0.01^2$ (m/s)$^2$ \\
        \hline
    \end{tabular}
    \label{tab: error sources}
\end{table}

The uncertainties associated with propulsive maneuver magnitude and direction and overall maneuver timing are based on typical operational errors resulting from thruster firing and attitude control. The absolute and relative navigation errors, $x_0$ and $\delta x_0$, are given element-wise in the 6D Cartesian frame (Earth-Centered-Inertial for absolute state, RTN for relative state) to match the general performance found in the literature \cite{giralo_digital}. To incorporate Cartesian values into (\ref{eq: dynamics}), a preliminary Gaussian Monte Carlo analysis is conducted. Each sample nonlinearly converts $\boldsymbol{\alpha}_c (t)$ or $a_c \delta \boldsymbol{\alpha}$ into the equivalent cartesian frame, generates and adds the applicable random Gaussian error, and converts new quantity back into its original frame. Again by law of large numbers, the covariance of these samples will converge to the true Keplerian orbital element and ROE covariance associated with navigation. Drag maneuver magnitude uncertainty is primarily caused by inaccurate atmospheric density modelling, which can vary by an order of magnitude in LEO. This analysis conservatively sets the 2$\sigma$ bound at $10^{-12}$kg/m$^3$, with the matching $\sigma$ value in Table \ref{tab: error sources}. Finally, the error in the direction of drag maneuvering is derived from the assumption that the chief and deputy share the chief's RTN frame. This error can exceed $8^\circ$ at along-track separations of 1000km. Therefore, the 2$\sigma$ bound is set conservatively again at $10^\circ$, resulting in a variance of $5^\circ$.

The covariance analysis for each error source applies the corresponding Gaussian distribution defined in Table \ref{tab: error sources} to the applicable term in the dynamics constraint given in (\ref{eq: linear err dynamics}) or (\ref{eq: nonlinear err dynamics}). For the linear errors, the final state covariance is given by solving (\ref{eq: linear error final state cov}). For the nonlinear errors, independent errors are sampled for each planned maneuver to find the resulting final state for each Monte Carlo sample. The final state covariance is the covariance of the final state samples. As mentioned before, the covariance analysis is conducted for each error source separately and only the resulting covariances are reported, as the linear transformation of a zero-mean Gaussian will always result in a zero-mean Gaussian. The performance bound for each reconfiguration can be found by solving (\ref{eq: performance bound}) with the corresponding final state covariance and a chi-squared cumulative probability of 95\%. The covariance analysis and worst-case performance bound results are shown in Tables \ref{tab: cov results} and \ref{tab: bound results}, where worst-case performance bounds are defined as the maximum bound for each dimension across the three reconfigurations.

\begin{table}[t]
    \centering
    \caption{Resulting maximum final state covariances from covariance analysis of hybrid maneuver plans solved for each IP dominance case candidate. Each error source is analyzed independently.}
    \begin{tabular}{c|c c c}
        $\epsilon$ & \multicolumn{3}{c}{max(diag($\boldsymbol{COV}_{\delta \boldsymbol{\alpha}_f}$))$^{1/2}$} [km] \\
        \hline \hline
       Dominance Case & $\delta a$ & $\delta \lambda$ & $\delta e$ \\
        \hline
       $|\boldsymbol{u}^{\text{prop}}_i|$ & $32.695$ & $49.943$ & $8.304$ \\
       $\measuredangle  \boldsymbol{u}^{\text{prop}}_i$ & $0.181$ & $3.833$ & $3.436$ \\
       $t_i^{\text{prop}}$ & $0.052$ & $0.023$ & $0.026$ \\
       $|\boldsymbol{u}^{\text{drag}}_i|$ & $1369.220$ & $685.046$ & $7.957$ \\
       $\measuredangle  \boldsymbol{u}^{\text{drag}}_j$ & $1.294$ & $0.668$ & $0.024$ \\
       $t_j^{\text{drag}}$ & $0.023$ & $0.034$ & $0.387$ \\
       $\boldsymbol{x}_0$ & $0.004$ & $0.009$ & $0.009$ \\
       $\delta \boldsymbol{x}_0$ & $\sim 10^{-7}$ & $\sim 10^{-7}$ & $\sim 10^{-7}$ \\
       \hline
    \end{tabular}
    \label{tab: cov results}
\end{table}

\begin{table}[b]
    \centering
    \caption{Resulting performance bounds from covariance analysis. Worst-case values are taken as largest bound for each dimension over all IP dominance case candidates. Each error source is analyzed independently.}
    \begin{tabular}{c|c c c c c c c c}
       \hspace{-2mm} $\epsilon$ \hspace{-2mm} & \multicolumn{8}{c}{Worst-Case ROE Performance Bounds $b_i$ [m]} \\
        \hline \hline
       \hspace{-2mm} $|\boldsymbol{u}^{\text{prop}}_i|$ \hspace{-2mm} & [ &\hspace{-5mm} 378 &\hspace{-3mm} 354,450 &\hspace{-3mm} 377 &\hspace{-3mm} 111 &\hspace{-3mm} 27 &\hspace{-3mm} 506 &\hspace{-5mm} ] \\
       \hspace{-2mm} $\measuredangle  \boldsymbol{u}^{\text{prop}}_i$ \hspace{-2mm} & [ &\hspace{-5mm} 42 &\hspace{-3mm} 27,203 &\hspace{-3mm} 12 &\hspace{-3mm} 53 &\hspace{-3mm} 23 &\hspace{-3mm} 44 &\hspace{-5mm} ] \\
       \hspace{-2mm} $t_i^{\text{prop}}$ \hspace{-2mm} & [ &\hspace{-5mm} $<$1 &\hspace{-3mm} 369 &\hspace{-3mm} 17 &\hspace{-3mm} 190 &\hspace{-3mm} 77 &\hspace{-3mm} 24 &\hspace{-5mm} ] \\
       \hspace{-2mm} $|\boldsymbol{u}^{\text{drag}}_i|$ \hspace{-2mm} & [ &\hspace{-5mm} 20,608 &\hspace{-3mm} 9,717,412 &\hspace{-3mm} 12,625 &\hspace{-3mm} 3,472 &\hspace{-3mm} $<$1 &\hspace{-3mm} 15,842 &\hspace{-5mm} ] \\
       \hspace{-2mm} $\measuredangle  \boldsymbol{u}^{\text{drag}}_j$ \hspace{-2mm} & [ &\hspace{-5mm} 19 &\hspace{-3mm} 9,184 &\hspace{-3mm} 34 &\hspace{-3mm} 55 &\hspace{-3mm} 41 &\hspace{-3mm} 25 &\hspace{-5mm} ] \\
       \hspace{-2mm} $t_j^{\text{drag}}$ \hspace{-2mm} & [ &\hspace{-5mm} 5 &\hspace{-3mm} 2,749 &\hspace{-3mm} 1 &\hspace{-3mm} 5 &\hspace{-3mm} $<$1 &\hspace{-3mm} 4 &\hspace{-5mm} ] \\
       \hspace{-2mm} $\boldsymbol{x}_0$ \hspace{-2mm} & [ &\hspace{-5mm} $<$1 &\hspace{-3mm} 61 &\hspace{-3mm} 10 &\hspace{-3mm} 64 &\hspace{-3mm} 31 &\hspace{-3mm} 9 &\hspace{-5mm} ] \\
       \hspace{-2mm} $\delta \boldsymbol{x}_0$ \hspace{-2mm} & [ &\hspace{-5mm} $<$1 &\hspace{-3mm} $<$1 &\hspace{-3mm} $<$1 &\hspace{-3mm} $<$1 &\hspace{-3mm} $<$1 &\hspace{-3mm} $<$1 &\hspace{-5mm} ] \\
       \hline
    \end{tabular}
    \label{tab: bound results}
\end{table}

Unlike the conclusions made by Chernick and D'Amico for propulsive-only control, the long control windows typical for hybrid control exacerbate the effect of nonlinear uncertainty sources on the final state covariance and confidence bound \cite{chernick}. With the exception of relative state navigation, all primary hybrid control errors considered here result in a significant deterioration of the final state accuracy, and general nonlinear error effects cannot be considered negligible with respect to linear error effects. However, the results do show that the largest contributor to final state error, by orders of magnitude, is the magnitude error of differential drag maneuvers, driven by uncertainty in the atmospheric density. Considering the differential drag magnitude errors in isolation, the bounds for Lyapunov stability given in (\ref{eq: converge bounds}) can be simplified to
\begin{equation} \label{eq: rho converge bounds}
    0 \leq \Bar{\rho}_c \leq 2\Tilde{\Bar{\rho}}_c,
\end{equation}
\noindent where $\Bar{\rho}_c$ is the average atmospheric density modelled by the generic reachable set-based solver (Alg. \ref{alg: architecture}, Lines 6 and 13) over the planned maneuver, and $\Tilde{\Bar{\rho}}_c$ is the actual time-integrated average of the density over a given drag maneuver. This constraint can be met by a simple tuning variable scalar within Lines 5 and 12 of Alg. \ref{alg: architecture}, posed as
\begin{equation}
    \Bar{\rho}_c = c_\rho \Bar{\rho}_{c,\text{model}},
\end{equation}
\noindent with coefficient $c_\rho$ on the time-averaged density from the on-board atmospheric density model $\Bar{\rho}_{c,\text{model}}$. Assuming that a lower bound can be placed on the value of atmospheric density at a given altitude, asymptotic Lyapunov stability can be guaranteed for the proposed re-solve architecture if the atmospheric density for differential drag maneuvers used by the solver (Alg. \ref{alg: architecture}, Lines 5 and 12) is no higher than double the lower bound. However, recall that Alg. \ref{alg: architecture} only re-solves the drag maneuver plan and requires drag to reject the other errors present in hybrid control. The most robust and simple strategy to ensure asymptotic Lyapunov convergence, while maintaining hybrid cost savings, is to model the effect of differential drag maneuvers using the lower bound for atmospheric density.

\section{Validation and Experimentation}

The validity of the proposed re-solve architecture and approach towards on-board modelled atmospheric density is demonstrated in this section for the SWARM-EX mission. SWARM-EX will study theromspheric and ionoshperic activity in LEO with a novel three 3U CubeSat swarm. Each identical spacecraft uses magnetorquers and reaction wheels for zero-propulsion 3-axis attitude determination and control and a single cold gas propulsion thruster for translational relative motion control. For the purposes of differential drag, the spacecraft are modelled to have a mass of 6kg, a coefficient of drag of 1.5, and a cross-sectional area between 0.01 - 0.09m$^2$, bounded by the dimensions of a 3U CubeSat solar panel array. These simulations will use two of the spacecraft in this swarm for simplicity, without loss of generality. The chief is assigned arbitrarily among the three identical spacecraft and is assumed to have the following initial mean orbital element state at the beginning of each simulation 
\begin{align}
    \boldsymbol{\alpha}_c &= 
    \begin{bmatrix}
        a & e & i & \Omega & \omega & M
    \end{bmatrix} \\
    &= 
    \begin{bmatrix}
        6798\text{km} & 0.003 & 51^\circ & 200^\circ & 70^\circ & 45^\circ
    \end{bmatrix}. \nonumber
\end{align}
\noindent The spacecraft are propagated by a 4th-order Runge-Kutta integration of the GVE over a 10s time step for propulsive-only and hybrid control, tailoring different perturbation models to the goals of each validation analysis. The propagation time step is increased to 40s for drag-only control, due to the significantly longer control windows required for drag in-plane reachability. Attitude constraints, such as Sun-pointing to remain power positive and Earth-pointing for atmospheric measurements, are neglected for simplicity, as they create a time dependency in the  $\Delta B^\rho$ bounds that do not affect the optimal drag maneuvers resulting from the reachable set theory analysis of the hybrid control problem in (\ref{eq: hybrid control problem}).

Two different reconfiguration scenarios will be considered for validation. In Reconfiguration 1, the spacecraft will reconfigure over a smaller inter-spacecraft separation, about 90km, over an 85 orbit control window, defined by the initial and final mean states
\begin{equation}
    \begin{bmatrix}
        a_c \delta \boldsymbol{\alpha}_0^1 \\ a_c \delta \boldsymbol{\alpha}_f^1
    \end{bmatrix} \hspace{-0.5mm} = \hspace{-0.5mm} 
    \begin{bmatrix}
        0.01 & \! \! \! \! -10 & \! \! 0.09 & \! \! 0.66 & \! \! 0.04 & \! \! 0.77 \\
        0 & \! \! \! \! -100 & \! \! 0 & \! \! 0.7 & \! \! 0 & \! \! 0.7
    \end{bmatrix} \text{km}.
\end{equation}
\noindent Reconfiguration 1 involves both a large change in $\delta \lambda$, or along-track separation, and $\delta \boldsymbol{e}$, returning to a nominal value for passively safe formation keeping. In Reconfiguration 2, the spacecraft will reconfigure over a large inter-spacecraft separation, about 500km, over a 200 orbit control window, defined by the initial and final mean states
\begin{equation}
    \hspace{-3mm}
    \begin{bmatrix}
        a_c \delta \boldsymbol{\alpha}_0^2 \\ a_c \delta \boldsymbol{\alpha}_f^2
    \end{bmatrix} \hspace{-0.5mm} = \hspace{-0.5mm} 
    \begin{bmatrix}
        0.125 & \! \! \! \! 800 & \! \! \! -0.05 & \! \! 0.1 & \! \! \! -0.01 & \! \! 0.08 \\
        0 & \! \! \! \! 1300 & \! \! \! 0 & \! \! 0 & \! \! \! 0 & \! \! 0
    \end{bmatrix} \text{km}. \hspace{-0.5mm}
\end{equation}
\noindent Reconfiguration 2 primarily involves a large change in $\delta \lambda$, minimally returning the remaining ROE to a nominal zero value. All control solutions (Alg. \ref{alg: architecture}, Lines 6 and 13) are generated by the reachable set theory-based hybrid control solver developed in previous work by the authors, allowing control solution computation, and the full proposed re-solve architecture, to run in linear time \cite{hunter_hybrid}. This is the same algorithm used in the previous section to conduct the hybrid control sensitivity covariance analysis. The architecture is designed to re-solve the drag maneuver plan every 5 orbits ($\Delta t_{\text{re-solve}} = $ 5 orbits), and deputies using the differential drag-only setup will solve their own control window on-board (Alg. \ref{alg: architecture}, Line 2).

\subsection{Experimental Bounds on Lyapunov Stability}

\begin{figure}
    \centering
    \includegraphics[scale=0.35]{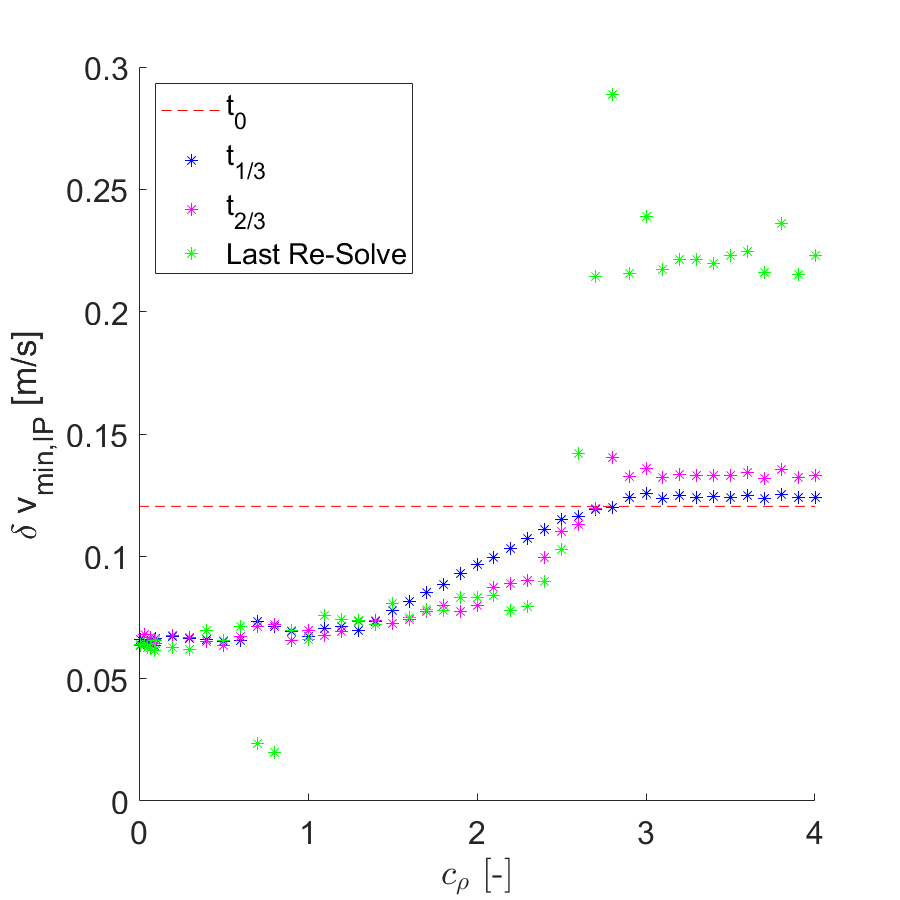}
    \caption{$\delta v_{\text{min,IP}}$ for Reconfiguration 1 at $t_{1/3}$, $t_{2/3}$, and the last re-solve of the control window, simulated under gravity and ideal drag. $c_\rho$ in Alg. 1 Lines 5 and 12 is varied to verify the derived bounds for Lyapunov stability. $\delta v_{\text{min,IP}}$ at $t_0$ is provided in red as the benchmark for convergence.}
    \label{fig: reconf1 converge bound validation}
\end{figure}

\begin{figure}
    \centering
    \includegraphics[scale=0.35]{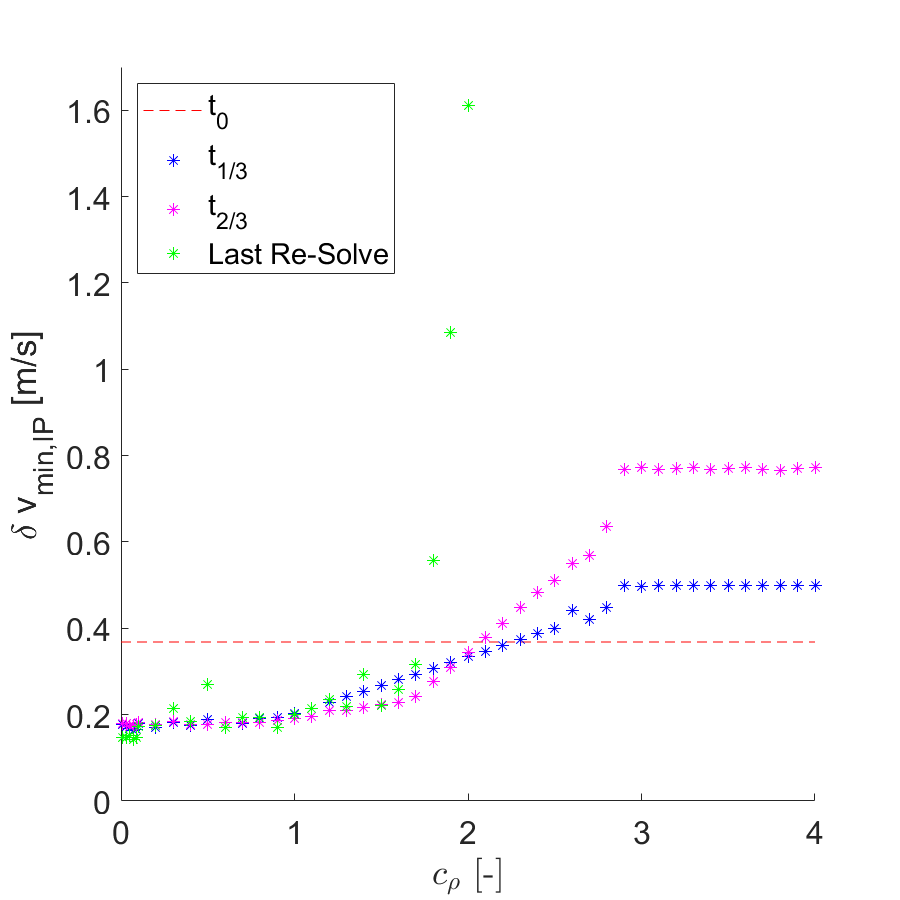}
    \caption{$\delta v_{\text{min,IP}}$ for Reconfiguration 2 at $t_{1/3}$, $t_{2/3}$, and the last re-solve of the control window, simulated under gravity and ideal drag. $c_\rho$ in Alg. 1 Lines 5 and 12 is varied to verify the derived bounds for Lyapunov stability. $\delta v_{\text{min,IP}}$ at $t_0$ is provided in red as the benchmark for convergence.}
    \label{fig: reconf2 converge bound validation}
\end{figure}

As a first step, the derived Lyapunov stability bounds on optimal reachable set theory-based control solutions (Alg. \ref{alg: architecture}, Lines 6 and 13) are verified in $S^3$, a high-fidelity orbit propagator designed by the Stanford Space Rendezvous Lab (SLAB) \cite{S3}. The validation reconfigurations are repeated with a variety of $c_\rho$ values (Alg. \ref{alg: architecture}, Lines 5 and 12), ranging from 0.01x to 4x the true mean value. To isolate the convergence of the designed choice of $c_\rho$, the simulation only incorporates perturbations from Earth's gravity (30x30 gravity model) and atmospheric drag using the NRLMSISE-00 density model \cite{nrlmsise00}. Ground truth atmospheric density is held constant, at the value at $t_0$, over the control window, and all spacecraft maneuver ideally with perfect on-board density knowledge ($\Bar{\rho}_{c,\text{model}} = \Tilde{\Bar{\rho}}$). Stable convergence is analyzed through the IP delta-v minimum $\delta v_{\text{min,IP}}$, a function of the current state relative to the final desired state and directly analogous to the Lyapunov function in (\ref{eq: V}). $\delta v_{\text{min,IP}}$ is found by solving (\ref{eq: pseudostate}) for the pseudostate at a given time in the control window, with the current state as $\delta \boldsymbol{\alpha}_0$ and the current time as $t_0$, and taking the maximum $\delta v_{\text{min}}$ from (\ref{eq: Chernick dv min da}), (\ref{eq: Chernick dv min dl}), and (\ref{eq: Chernick dv min de}) for the IP dominance case candidates. This analysis calculates the $\delta v_{\text{min,IP}}$ at three times in the control window to illustrate the evolution of $\delta v_{\text{min}}$ over time. The $\delta v_{\text{min,IP}}$ at the first third $t_{1/3}=t_0+1/3\tau$ and second third $t_{2/3}=t_0+2/3\tau$ of the control window duration $\tau=t_f-t_0$ will show the convergence properties of each $c_\rho$ throughout the control window while avoiding the time-based singularity in (\ref{eq: Chernick dv min dl}). The time of the last re-solve, approximately 5 orbits before the end of the control window, is used to demonstrate the effect of this singularity and verify convergence for the complete hybrid maneuver plan. The performance of the re-solve architecture is compared to the $\delta v_{\text{min,IP}}$ at the beginning of the control window $t_0$ where no control has yet occurred. Fig. \ref{fig: reconf1 converge bound validation} and Fig. \ref{fig: reconf2 converge bound validation} display the $\delta v_{\text{min,IP}}$ associated with each choice of modelled density.

Stable convergence for a value of $c_\rho$ in Alg. \ref{alg: architecture}, and the resulting hybrid control solution, is identified as those that decrease the $\delta v_{\text{min,IP}}$ from the value at $t_0$. In Reconfiguration 1, Alg. \ref{alg: architecture} achieves stable convergence for $c_\rho > 2$ at all three times considered, surpassing the derived bound of $\Bar{\rho}_c \leq 2\Tilde{\Bar{\rho}}_c$. For very low values of $c_\rho$, Alg. \ref{alg: architecture} stably converges to a residual non-zero $\delta v_{\text{min}}$ that represents the cost of the remaining IP propulsive maneuvers. These values are similar at all three times considered, because optimal IP maneuvers for large changes in $\delta \lambda$ occur in the first and last orbit of the control window, before and after all three considered times. The effect of the time-based singularity in (\ref{eq: Chernick dv min dl}) is evident for larger values of $c_\rho$, causing a spike in $\delta v_{\text{min,IP}}$ at the last re-solve. This singularity has a minimal effect on  smaller values of $c_\rho$ because the optimal cost of Reconfiguration 1 is dominated by the large change in $\delta \boldsymbol{e}$. In Reconfiguration 2, Alg. \ref{alg: architecture} achieves stable convergence for $c_\rho > 2$ at $t_{1/3}$ and $t_{2/3}$, but not at the last re-solve. This reconfiguration is primarily a large change in $\delta \lambda$, and the singularity effect in (\ref{eq: Chernick dv min dl}) can be seen over the entire control window, particularly at the last re-solve. This causes a large increase $\delta v_{\text{min}}$ at lower values of $c_\rho$ that skews the convergence properties of the full hybrid maneuver plan. Comparing these results to Reconfiguration 1, the time-based singularity effect of (\ref{eq: Chernick dv min dl}) primarily affects reconfigurations that are already $\delta \lambda$ dominant, and does not significantly change the output maneuver plan or convergence of the proposed re-solve architecture. Therefore, the experimental bounds for stable control in this analysis verify the derived Lyapunov stability bounds in (\ref{eq: rho converge bounds}). 

\subsection{Robust Performance}

The capabilities of the proposed re-solve architecture are demonstrated in $S^3$ with realistic operational and hardware errors. The full force simulation environment incorporates perturbations from Earth's gravity (30x30 gravity model), atmospheric drag using the NRLMSISE-00 density model, solar radiation pressure, and third body effects. The deputy solves and conducts control on-board with ideal maneuvering perturbed by the zero-mean Gaussian error models in Table \ref{tab: error sources}, with the exceptions of drag maneuver direction and magnitude. Ideal drag maneuver direction will be perturbed by inter-spacecraft separation, causing differences between the chief and deputy's RTN frames. Ideal drag maneuver magnitude will be  perturbed by variations between the on-board modelled time-average atmospheric density and the variable ground truth. At the inter-spacecraft separation in Reconfiguration 1, spacecraft absolute and relative states are assumed to be well characterized, and navigation uncertainties match the model in Table \ref{tab: error sources}. The spacecraft use the NRLMSISE-00 model as the reference truth on-board, the same model used as the simulation ground truth. At the inter-spacecraft separation in Reconfiguration 2, the increased separation is assumed to deteriorate the state estimation performance, and navigation uncertainties are taken from Table \ref{tab: error sources} and increased by an order of magnitude ($10\sigma$ or $100\sigma^2$). The spacecraft use the Harris-Priester atmospheric density model on-board, to simulate differences between the on-board reference truth density model and the ground truth \cite{harris-priester}. 

\begin{figure}[t]
    \centering
    \includegraphics[scale=0.15]{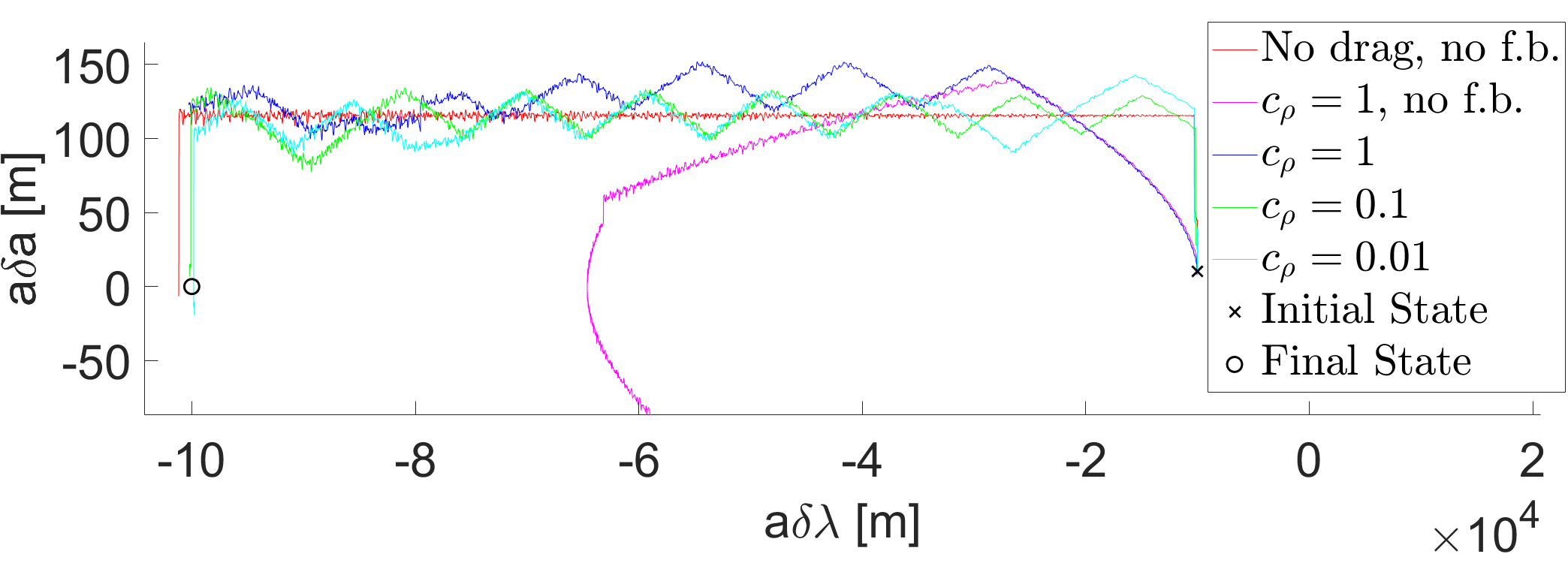}
    \caption{$\delta \boldsymbol{a}$ plane trajectory for Reconfiguration 1 for various control approaches and $c_\rho$ in Alg. 1 Lines 5 and 12 under full force dynamics model, operational errors, and nominal state estimation.}
    \label{fig: reconfig 1 hybrid}
\end{figure}

\begin{figure}[t]
    \centering
    \includegraphics[scale=0.15]{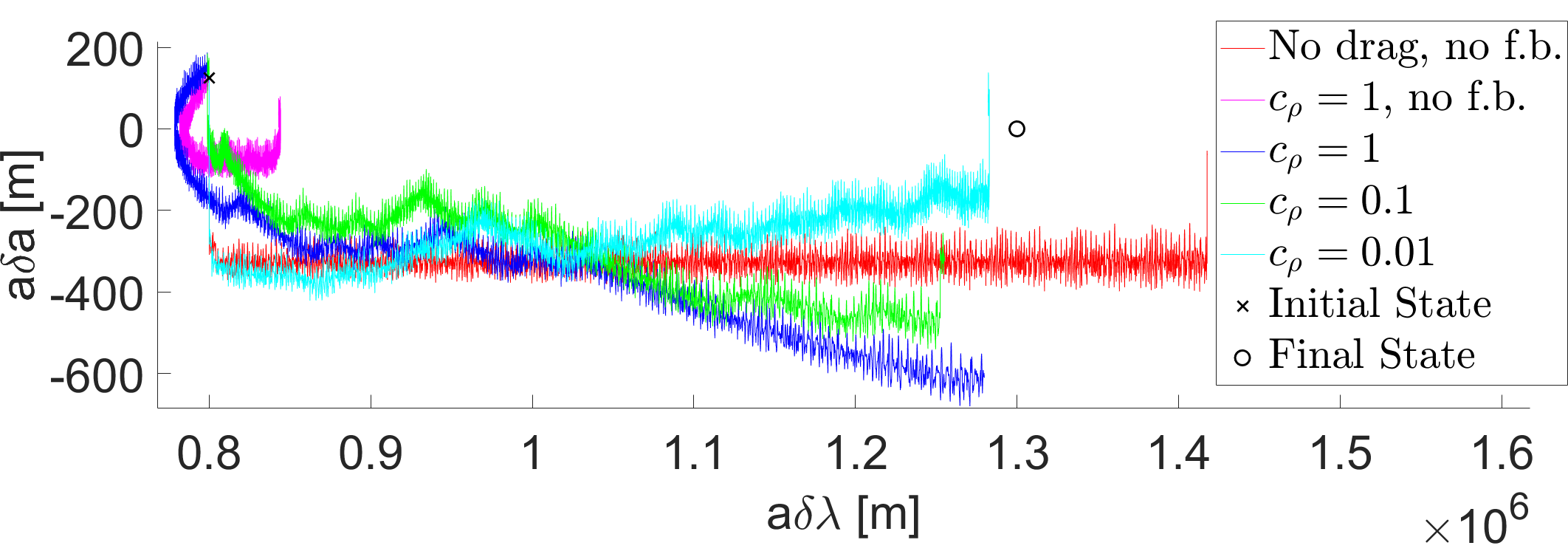}
    \caption{$\delta \boldsymbol{a}$ plane trajectory for Reconfiguration 2 for various control approaches and $c_\rho$ in Alg. 1 Lines 5 and 12 under full force dynamics model, operational errors, and poor state estimation.}
    \label{fig: reconfig 2 hybrid}
\end{figure}

The convergence and overall performance of the re-solve architecture is demonstrated comparatively against a single solve of the maneuver planning algorithm, conducted at the beginning of the control window without re-solves to provide a baseline to identify convergent control solutions. The single solve approach is equivalent to using Alg. \ref{alg: architecture} without Lines 7-16, the re-solve feedback loop, and is identified here as "no feedback" or "no f.b." For hybrid control, two "no feedback" baselines are used to evaluate robustness to uncertainty and dynamic modelling error. A single solve of the maneuver planning algorithm for propulsive-only ("no drag") control is compared against Alg. \ref{alg: architecture} hybrid maneuver planning to demonstrate the ability of the proposed architecture to achieve superior final state accuracy at a lower reconfiguration delta-v cost. A single solve of the maneuver planning algorithm for hybrid control and $c_\rho = 1$, compared against the single solve of propulsive-only control, shows the deterioration of final state accuracy caused by the inclusion of differential drag and poor characterization of the time-variable atmospheric density that Alg. \ref{alg: architecture} is designed to mitigate. The re-solve architecture is simulated with a variety of $c_\rho$ (Alg. \ref{alg: architecture}, Lines 5 and 12) to understand its behavior within the bounds of Lyapunov stable control. Tables \ref{tab: reconfig 1 hybrid} and \ref{tab: reconfig 2 hybrid} display the performance of each control approach for Reconfiguration 1 and 2, and Fig. \ref{fig: reconfig 1 hybrid} and \ref{fig: reconfig 2 hybrid} illustrate the trajectory of each control approach in the $\delta \boldsymbol{a}$ plane, where the effects of hybrid control and Alg. \ref{alg: architecture} are most evident.

\begin{table}[t]
    \centering
    \setlength{\tabcolsep}{3pt}
    \caption{IP cost and final ROE errors for Reconfiguration 1 for various control types and $c_\rho$ in Alg. 1 Lines 5 and 12 under full force dynamics model, operational errors, and nominal state estimation.}
    \begin{tabular}{c|c|c|c|c|c}
        Alg. \ref{alg: architecture} & IP $\delta v$ & \multicolumn{4}{c}{$a_c \delta \boldsymbol{\alpha}_f$ Error [m]} \\
        Drag Model & Cost [m/s] & $a_c\delta a$ & $a_c\delta \lambda$ & $a_c||\delta\boldsymbol{e}||_{2}$ & $a_c||\delta\boldsymbol{i}||_{2}$ \\
        \hline \hline
        \multicolumn{6}{l}{No Feedback} \\
        \hline
        No drag & 0.1226 & 6 & -2,138 & 6 & 12 \\
        $c_\rho = 1$ & 0.0125 & -86 & 41,143 & 19 & 74 \\
        \hline
        \multicolumn{6}{l}{Feedback} \\
        \hline
        $c_\rho = 1$ & 0.0123 & 115 & 453 & 120 & 5 \\
        $c_\rho = 0.1$ & 0.1074 & -4 & 561 & 7 & 15 \\
        $c_\rho = 0.01$ & 0.1223 & 3 & 635 & 3 & 23 \\
        \hline
    \end{tabular}
    \label{tab: reconfig 1 hybrid}
\end{table}

\begin{table}[t]
    \centering
    \setlength{\tabcolsep}{2pt}
    \caption{IP cost and final ROE errors for Reconfiguration 2 for various control types and $c_\rho$ in Alg. 1 Lines 5 and 12 under full force dynamics model, operational errors, and poor state estimation.}
    \vspace{-0.7mm}
    \begin{tabular}{c|c|c|c|c|c}
        Alg. \ref{alg: architecture} & IP $\delta v$ & \multicolumn{4}{c}{$a_c \delta \boldsymbol{\alpha}_f$ Error [m]} \\
        Drag Model & Cost [m/s] & $a_c\delta a$ & $a_c\delta \lambda$ & $a_c||\delta\boldsymbol{e}||_{2}$ & $a_c||\delta\boldsymbol{i}||_{2}$ \\
        \hline \hline
        \multicolumn{6}{l}{No Feedback} \\
        \hline
        No drag & 0.4114 & -65 & -117,758 & 56 & 285 \\
        $c_\rho = 1$ & 0.0142 & 6 & -457,010 & 79 & 673 \\
        \hline
        \multicolumn{6}{l}{Feedback} \\
        \hline
        $c_\rho = 1$ & 0.0041 & -610 & -20,181 & 117 & 113 \\
        $c_\rho = 0.1$ & 0.1926 & -324 & -44,798 & 196 & 179 \\
        $c_\rho = 0.01$ & 0.4073 & 65 & -17,814 & 90 & 24 \\
        \hline
    \end{tabular}
    \vspace{1.5mm}
    \label{tab: reconfig 2 hybrid}
\end{table}

The differences in operational errors, namely poor state estimation and atmospheric density prediction, are evident in the results for Reconfiguration 2, where final state errors are high even for low values of $c_\rho$. For both reconfigurations, the proposed re-solve architecture provides superior $\delta \lambda$ final state accuracy over both propulsive-only and hybrid control without feedback and achieves roughly equivalent or superior performance over the entire ROE final state for smaller values of $c_\rho$. Larger, stable values of $c_\rho$ do not guarantee decreased final state error across all ROE dimensions because the Lyapunov stability derived in this paper is proportional to the $\delta v_{\text{min}}$ in (\ref{eq: delta-v min}). Deputy trajectories in non-dominant dimensions may increase in tracking error as long as they do not increase the minimum cost of the reconfiguration. In cases where certain dimensions of a reconfiguration require a high level of precision (i.e. safe separation for collision avoidance), the desired final state accuracy may be achieved by decreasing $c_\rho$ significantly or restricting control to the plane of interest. Using a smaller $c_\rho$ in both Reconfiguration 1 and 2 improves performance in the $\delta a$ and $\delta \boldsymbol{e}$ final states but causes a corresponding increase in IP delta-v cost. This trade-off introduces a simple important tuning element for control design, adjusting the architecture to prioritize control accuracy or delta-v cost savings depending on the current needs of the mission. 

\begin{figure} [t]
    \centering
    \includegraphics[scale=0.15]{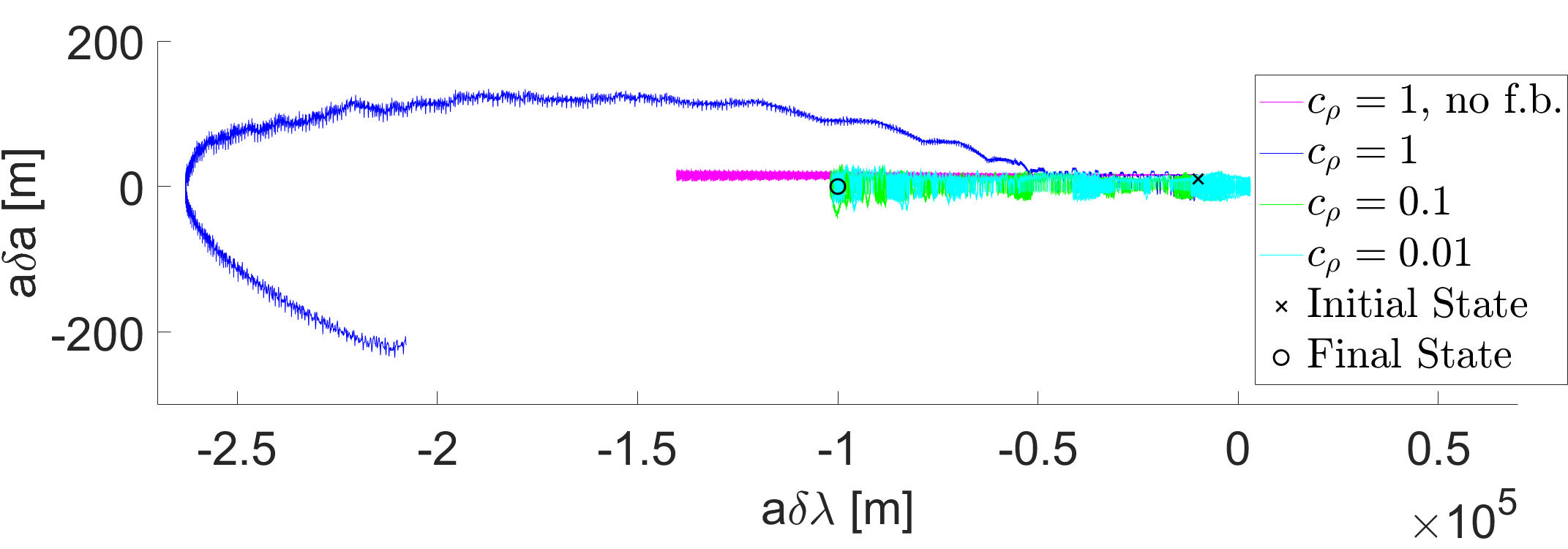}
    \caption{$\delta \boldsymbol{a}$ plane trajectory for Reconfiguration 1 for drag-only control and various $c_\rho$ in Alg. 1 Lines 5 and 12 under full force dynamics model, operational errors, and nominal state estimation.}
    \label{fig: reconfig 1 drag}
\end{figure}

\begin{figure} [t]
    \centering
    \includegraphics[scale=0.62]{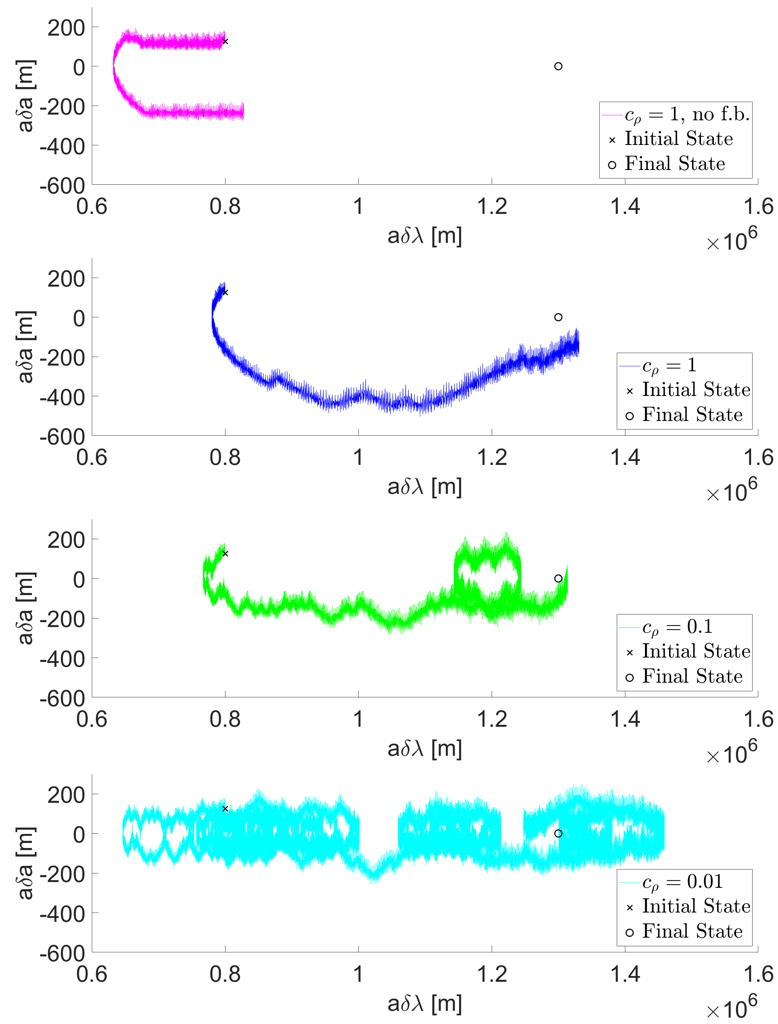}
    \caption{$\delta \boldsymbol{a}$ plane trajectory for Reconfiguration 2 for drag-only control and various $c_\rho$ in Alg. 1 Lines 5 and 12 under full force dynamics, operational errors, and poor state estimation. Each $c_\rho$ trajectory is split into separate plots for visual clarity.}
    \label{fig: reconfig 2 drag}
\end{figure}

For drag-only control, only the IP reconfiguration can be accomplished with differential drag control. The performance of Alg. \ref{alg: architecture} is evaluated comparatively for various values of $c_\rho$, measuring convergence to the final desired state and the control window required for modelled drag reachability. For Reconfiguration 1, $c_\rho = 0.01$ results in a control window $t_f-t_0$ of 111,614 hrs ($>$4,000 days) which is infeasible for a real mission and difficult to simulate numerically. To demonstrate the performance of small values of $c_\rho$, Reconfiguration 1 will use $c_\rho = 0.05$ instead of $c_\rho = 0.01$ with half of the control window given by (\ref{eq: drag dt}), the sufficient control window length to send all IP dominance case candidate $\delta v_{\text{min}}$ to zero without the buffer for navigation errors. A single solve of the maneuver planning algorithm for drag-only control is also included, to demonstrate the proposed re-solve architecture's rejection of dynamic modelling errors. Tables \ref{tab: reconfig 1 drag} and \ref{tab: reconfig 2 drag} detail the IP final state errors of Alg. \ref{alg: architecture} drag-only control for Reconfigurations 1 and 2, and Fig. \ref{fig: reconfig 1 drag} and Fig. \ref{fig: reconfig 2 drag} illustrate the trajectories in the $\delta \boldsymbol{a}$ plane, where convergence of drag-only control is again most evident.

The effect of increased estimation errors is again evident in the results for Reconfiguration 2, with high final state errors for low values of $c_\rho$. While Alg. \ref{alg: architecture} can provide superior final state accuracy over drag-only control without feedback, the re-solve architecture does not converge as reliably with drag-only control as with hybrid control. Drag maneuver magnitude is the primary source of error in hybrid control, and drag-only control exacerbates this effect, such that a larger reconfiguration in a particular ROE dimension may result in a larger final state error when feedback is not used. Reconfiguration 1 requires both a large $\delta \boldsymbol{e}$ and large $\delta \lambda$ reconfiguration, and Alg. \ref{alg: architecture} using $c_\rho = 1$ with feedback results in a lower $\delta \boldsymbol{e}$ final state error than Alg. \ref{alg: architecture} without feedback at the cost of higher $\delta a$ and $\delta \lambda$ errors. Again recall that the inclusion of feedback in Alg. \ref{alg: architecture} allows errors in non-dominant dimensions to increase so long as the reconfiguration optimal cost does not increase. Reconfiguration 2 primarily requires a large change in $\delta \lambda$, and the $\delta a$ and $\delta \boldsymbol{e}$ final state errors do not increase significantly when using feedback. Similar to hybrid control, a trade-off exists between overall final state accuracy and control window length, where a smaller $c_\rho$ results in a longer control window required for drag IP reachability. As demonstrated previously by the required modifications to drag-only validation for Reconfiguration 1, control window length quickly becomes infeasible for low values of $c_\rho$. The drag-only $c_\rho$ control tuning is therefore functionally limited at a value much higher than the lower convergence bound.

\begin{table}[t]
    \centering
    \setlength{\tabcolsep}{3pt}
    \caption{Control window length and final IP ROE errors for Reconfiguration 1 for drag-only control and various $c_\rho$ in Alg. \ref{alg: architecture} Lines 5 and 12 under full force dynamics model, operational errors, and nominal state estimation.}
    \begin{tabular}{c|c|c|c|c}
        Alg. \ref{alg: architecture} & Control Window & \multicolumn{3}{c}{$a_c \delta \boldsymbol{\alpha}_f$ Error [m]} \\
        Drag Model & $t_f-t_0$ [hrs] & $a_c\delta a$ & $a_c\delta \lambda$ & $a_c||\delta\boldsymbol{e}||_{2}$ \\
        \hline \hline
        \multicolumn{5}{l}{No Feedback} \\
        \hline
        $c_\rho = 1$ & 1,336 & 19 & -40,333 & 1,320  \\
        \hline
        \multicolumn{5}{l}{Feedback} \\
        \hline
        $c_\rho = 1$ & 1,336 & -217 & -107,792 & 776  \\
        $c_\rho = 0.1$ & 11,675 & 11 & 114 & 13  \\
        $c_\rho = 0.05$* & 11,445 & -17 & 848 & 11 \\
        \hline
        \multicolumn{5}{c}{*$t_f-t_0$ from (\ref{eq: drag dt}) halved for feasible simulation length}
    \end{tabular}
    \label{tab: reconfig 1 drag}
\end{table}

\section{Conclusions}

This paper presents a novel robust architecture for propulsive-differential drag "hybrid" control of Distributed Space Systems (DSS). The architecture leverages atmospheric drag, modulated by attitude at negligible propulsive cost, to both reduce the cost of satellite swarm reconfigurations and reject dynamic modelling errors that deteriorate final state accuracy. The algorithm begins with an initial solve of the full control problem and re-solves the remaining drag maneuver plan throughout the control window using the current state estimate. Through a Relative Orbital Element (ROE) parameterization, the optimal hybrid control problem is posed in a Linear Time Variant (LTV) formulation that is compatible with reachable set theory and generic reachable set-based solvers. Reachable set theory enables the optimal decoupling of the propulsive and differential drag control problems and provides a set of conditions for provable optimality that are extended in the design of the architecture to prove asymptotic Lyapunov stability for ideal control. The Lyapunov analysis also provides error bounds on stability for non-ideal control which can be enforced within the architecture on the inputs and outputs of the control solver.

\begin{table}[t]
    \centering
    \setlength{\tabcolsep}{2pt}
    \caption{Control window length and final IP ROE errors for Reconfiguration 2 for drag-only control and various $c_\rho$ in Alg. \ref{alg: architecture} Lines 5 and 12 under full force dynamics model, operational errors, and poor state estimation.}
    \begin{tabular}{c|c|c|c|c}
        Alg. \ref{alg: architecture} & Control Window & \multicolumn{3}{c}{$a_c \delta \boldsymbol{\alpha}_f$ Error [m]} \\
        Drag Model & $t_f-t_0$ [hrs] & $a_c\delta a$ & $a_c\delta \lambda$ & $a_c||\delta\boldsymbol{e}||_{2}$ \\
        \hline \hline
        \multicolumn{5}{l}{No Feedback} \\
        \hline
        $c_\rho = 1$ & 354 & 84 & -697,272 & 76  \\
        \hline
        \multicolumn{5}{l}{Feedback} \\
        \hline
        $c_\rho = 1$ & 369 & -137 & -31,101 & 45 \\
        $c_\rho = 0.1$ & 1,232 & -70 & 13,996 & 129 \\
        $c_\rho = 0.01$ & 10,151 & 27 & -52,421 & 65 \\
        \hline
    \end{tabular}
    \label{tab: reconfig 2 drag}
\end{table}

The primary operational error sources of hybrid control are detailed, and appropriate Gaussian error models are built for each source. A covariance sensitivity analysis is conducted to quantify the effect of each error source on the output maneuver plan and finds that the dominating source of error on hybrid control is the atmospheric density model. Therefore, the bounds on Lyapunov stable control can be simplified to bounds on the modelled time-averaged atmospheric density, such that the architecture asymptotically converges if the modelled effect of the maneuver is no greater than twice the actual effect of the maneuver. Asymptotic stability can be guaranteed if a lower bound can be placed on the value of atmospheric density at the altitude of the spacecraft swarm and is used as the atmospheric density value for control solutions.

The derived bounds on Lyapunov stability for reachable set theory-based optimal control are proven experimentally, applying the proposed re-solve architecture to an orbital environment perturbed by gravity and atmospheric drag. Full force simulations of the proposed re-solve architecture demonstrate its capability to significantly improve final state accuracy over no feedback maneuver planning approaches, while still maintaining the cost savings of differential drag. The performance of the architecture for both hybrid and drag-only control can be tuned by changing the atmospheric density used for control solutions. In hybrid control, decreasing the modelled density increases final state accuracy, but also increases IP cost. In drag-only control, decreasing the modelled density also increases final state accuracy, but increases the required control window and may require an infeasible control window length as the modelled density approaches the lower convergence bound.

Overall, this paper proves the asymptotic Lyapunov convergence of optimal control solvers derived from reachable set theory. This proof can be leveraged to create a simple re-solve strategy for robust cost-efficient control of DSS, designed by assessing the primary sources of dynamic error, analyzing each error's effect on final state accuracy, and conditioning the inputs and outputs of the reachable set-based solver. This approach is not specific to DSS and can be generalized to any generic LTV control problem with a norm-like cost function to achieve stable and provably optimal control despite an uncertain environment.

\acknowledgments
This work was supported by the National Science Foundation (Award Number 1936512). This project resulted from the Ideas Lab: Cross-cutting Initiative in CubeSat Innovations, an interdisciplinary program supported by Geosciences, Engineering, and Computer and Information Science and Engineering Directorates.


\bibliographystyle{IEEEtran}
\bibliography{references.bib}

\begin{biographywithpic}
{Simone D’Amico}{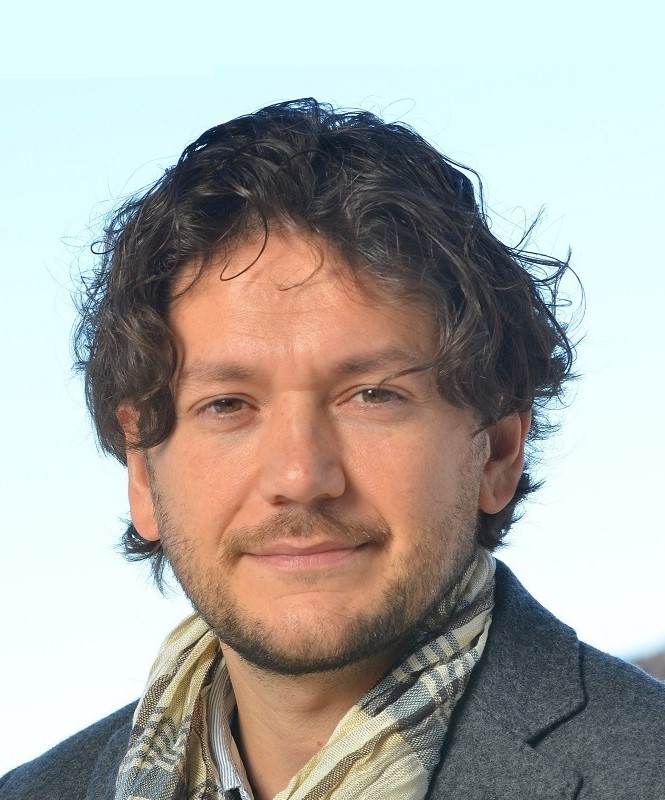}
receivedthe B.S. and M.S. degrees in aerospace engineering from Politecnico di Milano, Milan, Italy, in 2003, and the Ph.D. degree in Aerospace Engineering from the Delft University of Technology, Delft, The Netherlands, in 2010. He is currently an Associate Professor of Aeronautics and Astronautics (AA), W.M. Keck Faculty Scholar in the School of Engineering, and Professor of Geophysics (by Courtesy). He is the Founding Director of the Space Rendezvous Laboratory and Director of the AA Undergraduate Program. Before Stanford, Dr. D’Amico was research scientist and team leader at the German Aerospace Center (DLR) for 11 years. His research aims at enabling future miniature distributed space systems for unprecedented remote sensing, science, exploration and spaceflight sustainability. His efforts lie at the intersection of advanced astrodynamics, spacecraft navigation/control, and space system engineering. Dr. D’Amico is institutional PI of three upcoming autonomous satellite swarm missions funded by NASA and NSF. He is Fellow of AAS, Associate Fellow of AIAA, Associate Editor of AIAA JGCD, Advisor of NASA and space startups. He was the recipient of several awards, including Best Paper Awards at IAF (2022), IEEE (2021), AIAA (2021), AAS (2019) conferences, the Leonardo 500 Award by the Leonardo da Vinci Society/ISSNAF (2019), FAI/NAA’s Group Diploma of Honor (2018), DLR’s Sabbatical/Forschungssemester (2012) and Wissenschaft Preis (2006), and NASA’s Group Achievement Award for the GRACE mission (2004).
\end{biographywithpic}

\begin{biographywithpic}
{Matthew Hunter}{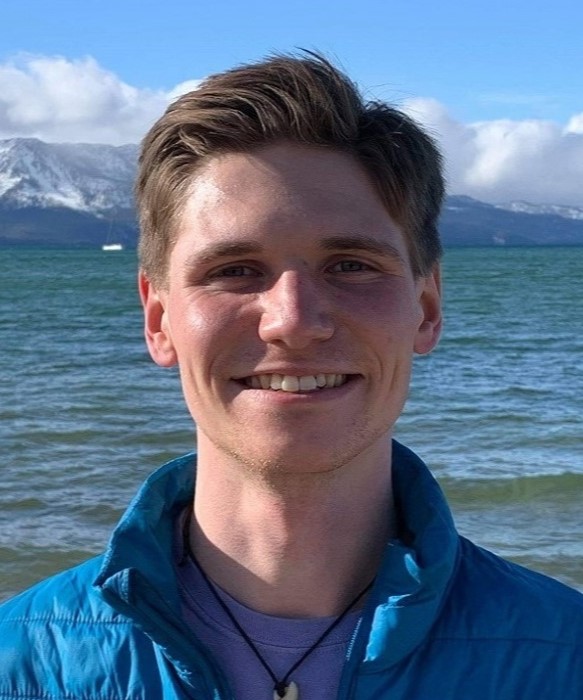}
is a Ph.D. student in the Space Rendezvous Laboratory. He “got out” of Georgia Tech with a Bachelor of Science in Mechanical Engineering, where he received the George W. Woodruff School of Mechanical Engineering Outstanding Scholar Award. He has interned at a variety of companies in the space and defense industry, including NASA Goddard, Raytheon Technologies, and Bell. Matthew is currently developing guidance and control algorithms for the Space Weather Atmospheric Multiscale EXperiment (SWARM-EX), and his research focuses on computationally efficient, optimal hybrid control techniques that incorporate perturbations to Keplerian orbit into propulsive maneuver schemes to decrease their delta-v cost. This work will enable fully integrated roto-translational spacecraft control to increase mission lifetime and enhance maneuver accuracy. 
\end{biographywithpic}

\end{document}